\documentclass[12pt]{article}
\newtheorem{theo}{{\bfseries Theorem}}[section]
\newtheorem{prop}[theo]{{\bfseries Proposition}}
\newtheorem{lem}[theo]{{\bfseries Lemma}}
\newtheorem{cor}[theo]{{\bfseries Corollary}}
\newtheorem{df}[theo]{{\bfseries Definition}}

\usepackage{amssymb,latexsym}
\usepackage[mathcal]{eucal}
\usepackage{amscd}
\usepackage{amsmath}
\usepackage{graphicx}
\usepackage{amsfonts}
\usepackage{amscd}

\numberwithin{equation}{section}

\begin{document}

\title{\bfseries  Dynamics of Discontinuous Maps Via Closed Relations}
\vspace{1cm}
\author{Ethan Akin\\
    Mathematics Department \\
    The City College \\
    137 Street and Convent Avenue \\
    New York City, NY 10031, USA}

    \vspace{.5cm}
\date{October, 2011}

\vspace{.5cm} \maketitle

\begin{abstract}  For the dynamics of a discontinuous map on a compact metric space, we describe an
approach using suitable closed relations and connect it with the continuous dynamics on an
invariant $G_{\delta}$ subset and with the continuous dynamics on the compact space of sample paths.\end{abstract}
\vspace{.5cm}


\vspace{.5cm}

This paper will appear in the Topology Proceedings.

\section{Introduction}

There has been some interest in extending the theory of dynamical systems to discontinuous maps with special focus
on quasi-continuous maps. See, for example,
Crannell and Martelli (2000). A subset $A$ of a metric space $X$ is called
 \emph{quasi-open} when $A \subset \overline{A^{\circ}}$. That is, $A$ is an open set together with part of its
 topological boundary. A map $f : X_1 \to X_2$ between metric spaces is call \emph{quasi-continuous} when
the pre-image of every open subset of $X_2$ is at least quasi-open.

As might be expected this leads to some oddities which, in my opinion, are best handled by using the existing theory
as extended to closed relations, e.g. as in Akin (1993), but with
a bit of trimming. Let me illustrate with an example of maps on
$I = [0,1]$.
\begin{equation}\label{1.1}
\begin{split}
f_0(x) \quad = \quad \begin{cases} \frac{1}{2} - x \qquad 0 \ \leq x \leq \frac{1}{2}, \\
\frac{3}{2} - x \qquad \frac{1}{2} < x \leq 1. \end{cases} \\ \\
f_1(x) \quad = \quad \begin{cases} \frac{1}{2} - x \qquad 0 \ \leq x < \frac{1}{2}, \\
\frac{3}{2} - x \qquad \frac{1}{2} \leq x \leq 1. \end{cases}
\end{split}
\end{equation}
These maps describe the flipping of each of the two intervals $[0,\frac{1}{2}]$ and $[\frac{1}{2},1]$ about its mid-point.
The maps disagree only in the choice of destination for the unique point of discontinuity $x = \frac{1}{2}$.

Iterating we obtain
\begin{equation}\label{1.2}
\begin{split}
f_0 \circ f_0(x) \quad = \quad f_0 \circ f_1(x) \quad = \quad \begin{cases} \ x \qquad 0 \ \leq x < 1, \\
\ 0 \ \qquad  x = 1. \end{cases} \\ \\
f_1\circ f_1(x) \quad = \quad f_1 \circ f_0(x)  \quad = \quad\begin{cases} \ 1 \ \qquad  \  x = 0, \\
 \ x \qquad 0 < x \leq 1. \end{cases}
\end{split}
\end{equation}

Instead, define
\begin{equation}\label{1.2a}
F_{01} \quad = \quad f_0 \cup f_1 \quad = \quad \overline{f_0} \quad = \quad \overline{f_1}.
\end{equation}
That is, $F_{01}$ is the common closure of $f_0$ and $f_1$, regarded as subsets of $I \times I$. If
we iterate the closed relation $F_{01}$ twice we obtain the identity map on all of $I$ together with the two extra points
$(0,1)$ and $(1,0)$.  That is,
\begin{equation}\label{1.3}
F_{01} \circ F_{01} \quad = \quad 1_I \ \cup \ \{ (0,1), (1,0) \}. \hspace{2cm}
\end{equation}

In my opinion the choice between $f_0$ and $f_1$ should be irrelevant to the dynamics and the way to handle this
is to use the closed relation $F_{01}$ instead. Next, the anomalous points should be discarded from
 our description of the dynamics of the system. There are three related ways to look at this approach and the
 purpose of this paper is to describe and relate them.

 A closed relation $F : X_1 \to X_2$ is a closed subset of the product $X_1 \times X_2$.  We regard a continuous map as a
 special case of a closed relation, identifying the map with what is sometimes called the graph of the map. Our spaces
 are compact metric spaces so that the composition of closed relations is a closed relation.  Thus, when $F$ is a closed
 relation on $X$, i.e. $X_1 = X_2 = X$, we can iterate defining inductively $F^{n+1} = F^n \circ F = F \circ F^n$.

 If $f : X_1 \to X_2$ is a continuous surjection of compact metric spaces then a closed subset $A \subset X_1$ is called
 \emph{minimal} for $f$ when it is minimal in the class of closed subsets of $X_1$ which are mapped onto $X_2$ by
 $f$.  Equivalently, $f(A) = X_2$ and if $B$ is a closed proper subset of $A$ then $f(B)$ is a proper subset of
 $X_2$. When $X_1$ itself is minimal for $f$ then the map is called \emph{irreducible}. This is equivalent to
 the condition that $f$ be an \emph{almost one-to-one} map, i.e. the set
 \begin{equation}\label{1.4}
 IN_{f} \quad = \quad \{ x \in X_1 : \ f^{-1}(f(x)) \ = \ \{ x \} \ \} \hspace{1cm}
 \end{equation}
 is dense in $X_1$.

 A continuous map $f : X_1 \to X_2$ is called  \emph{almost open} when
 for $A \subset X_1$, $A^{\circ} \not= \emptyset$ implies $f(A)^{\circ} \not= \emptyset$ where $A^{\circ}$
 is the interior of $A$.

 If $F : X_1 \to X_2$
  is a closed relation then by restricting the coordinate projections
 to $F$ we obtain the continuous maps
 $\pi_{1F} : F \to X_1$ and $\pi_{2F} : F \to X_2$.  We will call $F$  a \emph{suitable relation} when it satisfies
 two conditions:
 \begin{itemize}
 \item $\pi_{1F}$ is irreducible.
 \item $\pi_{2F}$ is almost open.
 \end{itemize}
 The first condition says that $F$ is a minimal subset for the projection $\pi_1 : X_1 \times X_2 \to X_1$.
 It says exactly that $F$ is the closure in $X_1 \times X_2$ of a continuous map
 defined on a dense subset of $X_1$. The second
 condition says that for $A \subset X_1$, $A^{\circ} \not= \emptyset$ implies $F(A)^{\circ} \not= \emptyset$.

With $F: X_1 \to X_2$ and $G : X_2 \to X_3$ suitable, it need not be true that the composition
$G \circ F : X_1 \to X_3$ is suitable.  In our example above, $F_{01}$ is a suitable relation but
$F_{01} \circ F_{01}$ is not suitable.

In general, if $F$ and $G$ are suitable then $G \circ F$ contains a unique subset minimal for $\pi_1$. We denote
this by $G \bullet F : X_1 \to X_3$. It clearly satisfies the first condition and in fact satisfies the second as well
and so is a suitable relation.  That is, $G \circ F$ contains a unique suitable relation from $X_1$ to $X_3$ and we
call this the \emph{suitable composition}.

The uniqueness requirement is why we impose the condition that $\pi_{2F}$ be almost open. Without it an open set might
be crushed to point and then under composition yield an open set of anomalous points.
Suppose we define $\tilde F_{01}$ on $[-1.1]$ by extending the definition
of $F_{01}$ by letting $x \mapsto \frac{1}{2}$ for all $x \in [-1,0]$. Then
\begin{equation}\label{1.5}
\tilde F_{01} \circ \tilde F_{01} \quad = \quad 1_I \ \cup \ [-1,0] \times \{ 0,1 \} \ \cup \ \{ (1,0) \}.
\end{equation}
This contains  many minimal subsets. If $U$ is any open subset of $[-1,0]$ then
$1_I \ \cup \ \overline{U} \times \{ 0 \} \ \cup \ ([-1,0] \setminus U) \times \{ 1 \}$ is a minimal subset
for the projection $\pi_1$.

Suitable composition is associative and so we obtain a category with objects compact metric spaces and morphisms
suitable relations with suitable composition. A suitable relation $F : X_1 \to X_2$ is an isomorphism in
this category exactly when $\pi_{2F}$ is an irreducible map.  In particular, an irreducible map is an isomorphism,
and any isomorphism is a composition of irreducible maps and the inverses of irreducible maps.

The relationship between quasi-continuous maps and closed relations has been described in an elegant paper by
Crannell, Frantz and LeMasurier (2006). From their results we see that:

\begin{theo}\label{theo1.1aa} Let $ X_1$ and $ X_2$  be compact metric spaces.
\begin{enumerate}
\item[(a)] If $g : X_1 \to X_2$ is a quasi-continuous function and $F$ is the closure of $g$ in $X_1 \times X_2$
then the closed relation $F$ is a minimal set for the projection $\pi_1 : X_1 \times X_2 \to X_1$, i.e. $\pi_{1F}$
is irreducible. Furthermore, $g$ is continuous at the points of a dense subset of $X_1$.

\item[(b)] If  $F \subset X_1 \times X_2$ is a closed relation
with $\pi_{1F} : F \to X_1$ irreducible and $g : X_1 \to X_2$ is a map with $g \subset F$,i.e.
\begin{equation}\label{1.4a}
g(x) \ \in \ F(x) \quad \mbox{ for all} \ \ x \in X_1,
\end{equation}
then $g$ is quasi-continuous and $F$ is the closure of $g$.
\end{enumerate}
\end{theo}

{\bfseries Remark:} A map $g : X_1 \to X_2$ which satisfies (\ref{1.4a}) is called a \emph{selection function} for $F$.
\vspace{.5cm}

{\bfseries Proof:} Clearly if $D$ is a dense subset of $X_1$ and a map $g_0 : D \to X_2$ is contained in
 a closed relation $F$ then $\overline{g_0}$ is a subset of $F$ which is mapped onto $X_1$ by
$\pi_1$. If $F$ is minimal for $\pi_1$ then $F = \overline{g_0}$. In particular, if $F$ is minimal then it
is the closure of any of its selection functions. Corollary 5 of Crannell, Frantz and LeMasurier (2006)
 says that these functions are quasi-continuous. This proves (b).

To say that a point $(x,y) \in F$ lies in $IN_{\pi_{1F}}$ says exactly that $F(x)$ is the singleton $\{ y \}$. In
particular, $g(x) = y$ for any selection function $g$. If $\pi_{1F}$ is irreducible then $IN_{\pi{1F}}$ is dense
in $F$ and so its projection $D$ is dense in $X_1$.  As we will show below (and it is easy to check) each
selection function $g$ is continuous at each point of this dense set $D$.

Now suppose that $g : X_1 \to X_2$ is quasi-continuous and that $F$ is its closure. Let $F_1$ be a closed subset
of $F$ which projects onto $X_1$ and let $h$ be a selection function for $F_1$.
Clearly, $\overline{h} \subset F_1 \subset F$. By Theorem 2 (ii) of Crannell, Frantz and LeMasurier (2006)
$\overline{h} = \overline{g}$. Hence, $F_1 = F$. Thus, $F$ is minimal for $\pi_1$. This completes the proof of (a).

$\Box$ \vspace{.5cm}

\begin{cor}\label{cor1.2aa} If $f : X_2 \to X_1$ is an irreducible map then any selection function
$g : X_1 \to X_2$ for the inverse relation $f^{-1} : X_1 \to X_2$ is a quasi-continuous map.\end{cor}

{\bfseries Proof:} The projection from the continuous map $f$ to its domain $X_2$ is a homeomorphism.
The projection of $f$ to its range
$X_1$ is thus homeomorphic to $f$ itself and so is irreducible. For the inverse relation $f^{-1}$ this means
that the projection to the domain $X_1$ is irreducible. By the above
theorem the selection functions are quasi-continuous.

$\Box$ \vspace{.5cm}

Closed relations are used in Crannell, Frantz and LeMasurier (2006) as a tool to study quasi-continuous functions.
Instead, I suggest that we focus on the closed relations themselves and eliminate the use of quasi-continuous
functions and the associated choices implicit in the use of selection functions.

Let $F$ be a suitable relation on $X$, i.e. $F : X \to X$. We will see that there exists a dense, $G_{\delta}$
subset $D^+_F$ of $X$ and a continuous map $t_F : D^+_F \to D^+_F$ such that $F$ is the
closure in $X \times X$ of the map $t_F$.
The three descriptions which I want to connect can be
labeled

\begin{itemize}
\item  Closing up the continuous map dynamics.

\item Suitable composition dynamics.

\item Sample path dynamics.
\end{itemize}

{\bfseries Closing up the continuous map dynamics:}  Since $D^+_F$ is a dense
$G_{\delta}$ subset of a compact metric space
it is a Polish space, i.e., a space that admits a complete, separable metric. Quite a bit of the usual dynamical systems
theory applies to $t_F$ on $D^+_F$ itself.  On the other hand, without compactness
there are undesirable holes in the theory.
For example, orbit sequences might have empty limit
point sets. Since $D^+_F$ is sitting in the compact space $X$ we can use the
limit points which occur in $X$. In general,
we focus on the dynamics of $t_F$ on $D^+_F$ and use the compactness of $X$ to define various limit point sets.

{\bfseries Suitable composition dynamics:}  We define the iterates
$F^{\bullet n+1} \ = \ F^{\bullet n} \bullet F \ = \  F \bullet F^{\bullet n}$.
That is,
we use the relation iterates $F^n$ and discard the points which do not lie in the unique minimal subset for
$\pi_{1F^n}$.

{\bfseries Sample path dynamics:} Let ${\mathbb N} = \{ 0,1,.... \}$.  For a closed relation $F$ on $X$ we define
the sample path space $X^+_F$ by
\begin{equation}\label{1.6}
X^+_F \quad = \quad \{ \ z \in X^{\mathbb N} : (z_i,z_{i+1}) \in F \quad \mbox{for all} \ \ i \in {\mathbb N} \ \}.
\end{equation}
This is a closed, subset of $X^{\mathbb N}$ which is $+$invariant
for the shift map $\sigma$ on $X^{\mathbb N}$ defined by
\begin{equation}\label{1.7}
\sigma(z)_i \quad = \quad z_{i+1} \qquad \mbox{for all} \ \ i \in {\mathbb N}. \hspace{2cm}
\end{equation}
The map $\pi_0 \times \pi_n$ maps $X^+_F$ onto the $n$-fold iterate $F^n$. Equivalently, $\pi_0$ maps
the restriction of the map $\sigma^n|X^+_F$ onto the $n$-fold relation composite $F^n$.

When $F$ is a suitable relation, there is a unique subset of $X^+_F$ which is minimal for $\pi_0$. We label it $S^+_F$
Thus, the restriction $\pi_0 : S^+_F \to X$ is irreducible.
This subset is $+$invariant for the shift and
 we denote by $s_F : S^+_F \to S^+_F$ the restriction of the shift map. It is an almost open, continuous map.

 We will show that that three approaches are related in that the closure in $X \times X$ of the $n$-fold iterate
 $(t_F)^n \subset D^+_F \times D^+_F$ is exactly $F^{\bullet n}$.  Furthermore, the map $\pi_0 \times \pi_n$ maps $S^+_F$
 onto $F^{\bullet n}$. Equivalently, $\pi_0$ maps
 $(s_F)^n$ on $S^+_F$ onto  $F^{\bullet n}$.

 When $F$ is a suitable relations isomorphism, there is an irreducible lift to a homeomorphism.

 Finally, recall that all measure spaces associated with nonatomic measures
 are isomorphic to the unit interval with Lebesgue measure.
 Similarly, every compact metric space without isolated points is isomorphic
 in the suitable relations category to the Cantor Set.
 Furthermore any suitable relation on such a space can be lifted, via an
 irreducible map to an almost open, continuous map on the Cantor Set and
 a suitable relations isomorphism can be similarly lifted to a homeomorphism on the Cantor Set.
 \vspace{1cm}

 \section{Almost Open Maps and Closed Relations}

 All of our spaces are nonempty separable metric spaces. For a subset $A$ we will use the usual
 notation $\overline{A}$ and $A^{\circ}$ for the closure
 and interior, respectively. Our main focus is upon compact metric spaces.
 We will reserve the letters $X, Y$ with or without subscripts for such
 spaces and use other letters for more general metric spaces.

 A  space is Polish when it is separable and admits a complete metric.
 For example, a compact metric space is Polish.

 In this section we review some foundational results about almost open maps and about closed relations between compacta.

 \begin{df}\label{df2.1}  Let $f : E_1 \to E_2$ be a continuous map.
 \begin{itemize}
 \item The map $f$ is \emph{open at} $x \in E_1$ if $x \in U$ implies $f(x) \in f(U)^{\circ}$ for all open
 $U \subset E_1$.

 \item The map $f$ is \emph{open} if $f(U)$ is open for all open $U \subset E_1$.

 \item The map $f$ is \emph{almost open} if $f(U)^{\circ} \not= \emptyset$ for all
 nonempty open $U \subset E_1$.

 \item  The map $f$ is \emph{weakly open} if $( \ \overline{f(U)} \ )^{\circ} \not= \emptyset$ for all
 nonempty open  $U \subset E_1$.
 \end{itemize}
 \end{df}
 \vspace{.5cm}

 Clearly, open implies almost open and almost open implies weakly open.  We label by
 $OPEN_f$ the set of points of $E_1$ at which $f$ is open. For an open  $U \subset E_1$ we define the
 open $U_f \subset U$ by
 \begin{equation}\label{2.1}
 U_f \quad =_{def} \quad U \cap f^{-1}[ (f(U)^{\circ} ]. \hspace{2cm}
 \end{equation}
 That is, $U_f$ is the set of points of $x \in U$ such that $f(U)$ is a neighborhood of $f(x)$.

 \begin{prop}\label{prop2.2}  Let $f : E_1 \to E_2$ be a continuous map.
\begin{enumerate}
\item[(a)] The map $f$ is open at $x$ iff $x \in U_f$ for every open $U \subset E_1$ with $x \in U$.

\item[(b)] The map $f$ is open iff $U = U_f$ for every open $U \subset E_1$.

\item[(c)] The map $f$ is open iff it is open at every point of $E_1$, i.e. $OPEN_f = E_1$.

\item[(d)] The map $f$ is almost open iff $U_f$ is dense $U$ for every open $U \subset E_1$.

\item[(e)] The map $f$ is almost open if it is open at a dense set of points, i.e. if $OPEN_f$ is dense in $E_1$.
\end{enumerate}
\end{prop}

{\bfseries Proof:} Results (a), (b) and (c) are obvious.  It is also clear that $f$ is almost open iff
$U_f \not= \emptyset$ whenever $U$ is open and nonempty. This, together with (a) implies (e).
Density of $U_f$ in $U$ certainly implies that $U_f$ is nonempty when $U$ is.

We are left with showing that $U_f$ is dense in $U$ whenever $f$ is almost open. Let
$V$ be an arbitrary nonempty open subset of $U$.  Because $f$ is almost open, $V_f \subset U_f \cap V$ is nonempty.
Hence, there exist points of $U_f$ which lie in $V$.

$\Box$ \vspace{.5cm}

 \begin{prop}\label{prop2.2a}  Let $f : E_1 \to E_2$ be a continuous map.

 The map $f$ is almost open iff $D $ dense in $E_2$ implies $f^{-1}(D)$ is dense in $E_1$.

 The map $f$ is weakly open iff $D $ open and dense in $E_2$ implies $f^{-1}(D)$ is open and dense in $E_1$.

 If $E_1$ is a Polish space and $f$ is weakly open then $D$ a dense $G_{\delta}$ subset of $E_2$ implies $f^{-1}(D)$
 is a dense $G_{\delta}$  subset of $E_1$.
 \end{prop}

 {\bfseries Proof:}  If $f$ is not almost open then there exists a nonempty open $U$ with $f(U)^{\circ}$ empty
 and so the complement $D$ of $f(U)$ in $ E_2$  is dense with $f^{-1}(D)$ disjoint from $U$. If $f$ is not weakly open
 then there exists $U$ with $( \ \overline{f(U)} \ )^{\circ}$ empty and so the complement $D$ of $\overline{f(U)}$
 is open and dense with $f^{-1}(D)$ disjoint from $U$.

 For the converses, let $D$ be a dense subset of $E_2$ and $U$ be an arbitrary nonempty open subset of $E_1$.
 If $f$ is almost open then $f(U)$ has a nonempty interior and so meets $D$.  Hence, $U$ meets $f^{-1}(D)$.
 If $f$ is weakly open then $\overline{f(U)}$ has a nonempty interior and so meets $D$.  If $D$ is also open then
 $f(U)$ meets $D$ and so, as before, $U$ meets $f^{-1}(D)$.

 If $D$ is a countable intersection of a sequence  $\{ V_n \}$ of open subsets of $E_2$ then each $V_n$ is a dense
 open set and so when $f$ is weakly open each $f^{-1}(V_n)$ is a dense open subset of $E_2$.  Thus,
 $f^{-1}(D) = \bigcap_n \ f^{-1}(V_n)$ is a $G_{\delta}$ which is dense by the Baire Category Theorem
when $E_1$ is Polish.

 $\Box$ \vspace{.5cm}

 \begin{cor}\label{cor2.2b} Let $f : E_1 \to E_2$ and $g : E_2 \to X_3$ be continuous maps.

 If both $f$ and $g$ are almost open (or weakly open) then $g \circ f$ is almost open (resp. weakly open).

 If $f$ is surjective and $g \circ f$ is almost open (or weakly open) then $g$ is almost open (resp. weakly open).
 \end{cor}

 {\bfseries Proof:}  Let $D$ be dense in $X_3$.  If $f$ and $g$ are almost open then $g^{-1}(D)$ is dense in
 $E_2$ and so $f^{-1}(g^{-1}(D)) = (g \circ f)^{-1}(D)$ is dense in $E_1$.  So $g \circ f$ is almost open.

 If $g \circ f$ is almost open then $f^{-1}(g^{-1}(D)) = (g \circ f)^{-1}(D)$ is dense in $E_1$ and since
 $f$ is surjective $g^{-1}(D) = f(f^{-1}(g^{-1}(D)))$ is dense in $E_2$.  Hence $g$ is almost open.

 For weakly open the proof is the same with dense replaced by open and dense throughout.

 $\Box$ \vspace{.5cm}

\begin{theo}\label{theo2.3} Let $f : E_1 \to E_2$ be a continuous map with $E_1$
locally compact.

The map $f$ is almost open iff it is weakly open.
\end{theo}

{\bfseries Proof:}  If $f$ is weakly open, $U$ is open and $x \in U$ then there exists a compact
$K \subset U$ with $x \in K^{\circ}$. Since $f(K)$ is compact, it is closed.
Because $f$ is weakly open $[\overline{f(K^{\circ})}]^{\circ}$ is nonempty.  It is contained in
$f(K)^{\circ} \subset f(U)^{\circ}$. As $U$ was arbitrary, it follows that $f$ is almost open.

$\Box$ \vspace{.5cm}

While our focus will be on almost open maps, the weaker notion is a useful tool
because of the following result which we will use often
and refer to as the \emph{Variation of Domain  Theorem}

\begin{theo}\label{theo2.4}   Let $f : E_1 \to E_2$ be a continuous map
and for $i = 1,2$ let $D_i$ be a
dense subset of $E_i$.  Assume that $f(D_1) \subset D_2$ so that we can define the restricted map
$f : D_1 \to D_2$.

The map $f : E_1 \to E_2$ is weakly open iff the restriction $f : D_1 \to D_2$ is weakly open.
\end{theo}

The proof requires the following

\begin{lem}\label{lem2.5} Let $D$ be a dense subset of $E$.
\begin{enumerate}
\item[(a)] If $U \subset E$ is open then $U \cap D$ is dense in $U$.

\item[(b)] If $A \subset D$ then $\overline{A} \cap D$ is the $D$-closure of $A$.

\item[(c)] If $B \subset E$ is closed then $B^{\circ} \cap D$ is the $D$-interior of $B \cap D$.
\end{enumerate}
\end{lem}

{\bfseries Proof:} (a) If $V \subset U$ is an arbitrary nonempty open subset then $V$  meets
the dense set $D$.

(b) $\overline{A} \cap D$ is $D$-closed and contains $A$ and so contains the $D$-closure.  On the other hand, there
is a closed set $F$ such that $F \cap D$ is the $D$-closure of $A$. In particular, $F$ contains $A$ and so
since $F$ is closed it contains $\overline{A}$. Hence the $D$-closure of $A$ which is $F \cap D$ contains
$\overline{A} \cap D$.

(c) $B^{\circ} \cap D$ is $D$-open and is contained in $B \cap D$ and so is contained in the $D$-interior of
$B \cap D$. On the other hand, there exists an open set $U$ such that $U \cap D$ is the $D$-interior of
$B \cap D$. Since $U \cap D \subset B$ and $B$ is closed, $\overline{U \cap D} \subset B$. By (a) $U \cap D$ is
dense in $U$ and so $\overline{U} = \overline{U \cap D} \subset B$. Hence, the open set $U$ is contained in $B$ and
so in
$B^{\circ}$.  Hence the $D$-interior of $B \cap D$ which is $U \cap D$  is contained in $B^{\circ} \cap D$.

$\Box$ \vspace{.5cm}

{\bfseries Proof of Theorem \ref{theo2.4}:}
 The open subsets of $D_1$ are of the form $U \cap D_1$.  By Lemma \ref{lem2.5}(a) $U \cap D_1$
is dense in $U$ and so $f(U \cap D_1)$ is dense in $f(U)$. Thus,
\begin{equation}\label{2.2}
\overline{f(U)} \quad = \quad \overline{f(U \cap D_1)} \hspace{3cm}
\end{equation}

Since $f(U \cap D_1) \subset D_2$, Lemma \ref{lem2.5}(b) implies that
$\overline{f(U)} \cap D_2 = \overline{f(U \cap D_1)} \cap D_2$ is the $D_2$-closure of $f(U \cap D_1)$.

Finally Lemma \ref{lem2.5} (c) implies that $[\overline{f(U)}]^{\circ} \cap D_2$ is the $D_2$-interior
of $\overline{f(U \cap D_1)} \cap D_2$ which we have seen is the $D_2$-closure of $f(U \cap D_1)$.

Thus,  $[\overline{f(U)}]^{\circ} \cap D_2$ is the $D_2$-interior of the $D_2$-closure of $f(U \cap D_1)$.

By Lemma \ref{lem2.5}(a) again $[\overline{f(U)}]^{\circ} \cap D_2$ is dense in $[\overline{f(U)}]^{\circ}$.

It follows that$[\overline{f(U)}]^{\circ} $ is nonempty   iff
the $D_2$-interior of the $D_2$ closure of $f(U \cap D_1)$ is nonempty. As $U$ varies over all the nonempty open
subsets of $E_1$, $U \cap D_1$ varies over all the nonempty $D_1$-open subsets of $D_1$. Thus, $f$ is
weakly open iff the restriction $f : D_1 \to D_2$ is weakly open.

$\Box$ \vspace{.5cm}

{\bfseries Remark:}  Note that the identity map on the unit interval $I$ is open, but the restriction with
$D_1$ the rationals and $D_2 = I$ is not almost open.
\vspace{.5cm}

\begin{prop}\label{prop2.6} Let $f : E_1 \to E_2$ be a continuous map with $E_1$ a Polish space.

The set $OPEN_f$, the set of points of $E_1$ at which $f$ is open, is a $G_{\delta}$ subset of $E_1$. The map $f$ is almost open iff $OPEN_f$ is dense in $E_1$.
\end{prop}

{\bfseries Proof:} Let $\mathcal{A}$ be the set of pairs $(U_1,U_2)$ members of a countable basis
with $\overline{U_2} \subset U_1$ so that $\{ U_1, E_1 \setminus \overline{U_2} \}$ is an open cover of $E_1$.
I claim that
\begin{equation}\label{2.3}
OPEN_f \quad = \quad \bigcap_{(U_1,U_2) \in \mathcal{A}} \ (U_1)_f \cup (E_1 \setminus \overline{U_2})_f. \hspace{2cm}
\end{equation}
By Proposition \ref{prop2.2}(a) $x$ is in the intersection when $f$ is open at $x$.  On the other hand, suppose that
$x$ is in the above intersection and $U$ is a open set with $x \in U$. We can choose $(U_1,U_2) \in \mathcal{A}$
such that $x \in U_2 $ and $U_1 \subset U$. Since $x$ is in the intersection but not in $E_1 \setminus \overline{U_2}$
it follows that $x \in (U_1)_f \subset U_f$. As $U$ was arbitrary, $f$ is open at $x$.

Now if $f$ is almost open then by Proposition \ref{prop2.2}(d) each of the open sets
$(U_1)_f \cup (E_1 \setminus \overline{U_2})_f$ is dense in $(U_1) \cup (E_1 \setminus \overline{U_2}) = E_1$.
As the countable intersection of dense open sets, $OPEN_f$ is dense by the
Baire Category Theorem applied to the Polish space $E_1$.

The converse follows from Proposition \ref{prop2.2}(e).

$\Box$ \vspace{.5cm}

We now review the theory of relations following Akin (1993).

As do the set theorists, we regard a map $f : E_1 \to E_2$ as a subset of $E_1 \times E_2$. For example,
the identity map $1_E $ is the diagonal set $ \{ (x,x): x \in E \}$. A relation $F : E_1 \to E_2$ is an arbitrary
subset of $E_1 \times E_2$. For $x \in E_1$  and $A \subset E_1$ we let
\begin{equation}\label{2.4}
\begin{split}
F(x) \quad = \quad \{ y \in E_2 : (x,y) \in F \} \hspace{3cm} \\
F(A) \quad = \quad \{ y : (x,y) \in F \quad \mbox{for some} \ x \in A \} \hspace{.5cm}  \\
= \quad \bigcup_{a \in A} \ F(a) \quad = \quad \pi_1(F \cap (A \times X_2)), \hspace{1cm}
\end{split}
\end{equation}
where $\pi_1 : X_1 \times X_2  \to X_1$ is the projection to the first coordinate.

Thus, $F$ is a map exactly when each $F(x)$ is a singleton.  In that case, we use $F(x)$ to stand both for the
singleton set and for the point it contains, allowing context to determine the choice of meaning.

The \emph{inverse relation} $F^{-1}: E_2 \to E_1$ is defined by
\begin{equation}\label{2.5}
F^{-1} \quad = \quad \{ (y,x) : (x,y) \in F \}. \hspace{3cm}
\end{equation}
Thus, for $B \subset E_2$ we see that
\begin{equation}\label{2.6}
F^{-1}(B) \quad = \quad \{ x \in E_1 : F(x) \cap B \not= \emptyset \}. \hspace{2cm}
\end{equation}

We call
\begin{equation}\label{2.7}
F^{-1}(E_2) \quad = \quad \{ x \in E_1 : F(x)  \not= \emptyset \} \hspace{2cm}
\end{equation}
the \emph{domain} of $F$. When $F^{-1}(E_2) = E_1$ we will say that \emph{$F$ has full domain} or
call $F$ a \emph{full domain relation}.

For $A_1 \subset E_1$ and $A_2 \subset E_2$ the \emph{restriction} is $F \cap (A_1 \times A_2)$ regarded as a relation
from $A_1$ to $A_2$. Thus, we can always restrict to obtain a full domain relation.

There is an alternative version of the pre-inverse of $B$ which we will label $F^*(B)$.
\begin{equation}\label{2.8}
\begin{split}
F^{*}(B) \quad = \quad \{ x \in E_1 : F(x) \subset B  \}  \hspace{2cm} \\
\quad = \quad E_1 \setminus F^{-1}(E_2 \setminus B). \hspace{3cm}
\end{split}
\end{equation}

Both operators are monotone with respect to set inclusion.  While $F^{-1}(\emptyset) = \emptyset$,
$F^*(\emptyset)$ is the complement of the domain of $F$, i.e. the set of $x$ such that $F(x) = \emptyset$.
Clearly,
\begin{equation}\label{2.9}
F^*(B) \quad \subset \quad F^{-1}(B) \cup F^*(\emptyset).  \hspace{2cm}
\end{equation}
So $F^*(B) \subset F^{-1}(B)$ for all $B$ exactly when $F$ is a full domain relation.

If $\{ B_{\alpha} \}$ is a family of subsets of $E_2$ then
\begin{equation}\label{2.10}
\begin{split}
F^{-1}(\bigcup_{\alpha} \ B_{\alpha}) \quad = \quad \bigcup_{\alpha} \ F^{-1}( B_{\alpha}), \hspace{1cm}\\
F^{*}(\bigcap_{\alpha} \ B_{\alpha}) \quad = \quad \bigcap_{\alpha} \ F^{*}( B_{\alpha}), \hspace{1cm}.
\end{split}
\end{equation}

When $F$ is a function these two operators agree and are the usual pre-image operator.

If $F : E_1 \to E_2$ and $G : E_2 \to E_3$ then we define
\begin{equation}\label{2.11}
\begin{split}
G \otimes F \quad =_{def}\quad (F \times E_3) \cap (E_1 \times G) \quad = \hspace{2cm}\\
 \{ (x,y,z) \in E_1 \times E_2 \times E_3 : (x,y) \in F \ \mbox{and} \ (y,z) \in G \}.
\end{split}
\end{equation}
and the \emph{composition} $G \circ F$
\begin{equation}\label{2.12}
\begin{split}
G \circ F \quad =_{def}\quad \pi_{13}(G \otimes F) \quad = \quad \hspace{3cm}\\
 \{ (x,z) \in E_1  \times E_3 : \mbox{for some} \ y \in E_2 \  \ (x,y) \in F \ \mbox{and} \ (y,z) \in G \}
\end{split}
\end{equation}
where $\pi_{13}$ is the projection to the product of the first and third coordinates.

Clearly, for $A \subset E_1$
\begin{equation}\label{2.13}
(G \circ F)(A) \quad = \quad G(F(A))
\end{equation}
and
\begin{equation}\label{2.14}
(G \circ F)^{-1} \quad = \quad F^{-1} \circ G^{-1}.
\end{equation}
Hence, for $B \subset E_3$
\begin{equation}\label{2.15}
\begin{split}
(G \circ F)^{-1}(B) \quad = \quad F^{-1}(G^{-1}(B)),\hspace{2cm} \\
(G \circ F)^{*}(B) \quad = \quad F^{*}(G^{*}(B)),\hspace{2.5cm}
\end{split}
\end{equation}

Composition of relations automatically takes care of domain problems.
Suppose $D_1 \subset E_1$, $D_2 \subset E_2$ on which are
defined maps $f : D_1 \to E_2$ and $g : D_2 \to E_3$. We can restrict
$f$ to $f^{-1}(D_2) \subset D_1$ to get a map $f : f^{-1}(D_2) \to D_2$
which we can compose with $g$ to get a map $g \circ f : f^{-1}(D_2) \to E_3$. Instead we can regard
$f : E_1 \to E_2$ and $g : E_2 \to E_3$ as relations
with $f^{-1}(E_2) = D_1$ and $g^{-1}(E_3) = D_2$. The composed relation
$g \circ f : E_1 \to E_3$ has $(g \circ f)^{-1}(E_3) = f^{-1}(D_2)$ and is
just the composed function on that set.

Since a relation $F : E_1 \to E_2$ is a subset of the product, we can define
the maps $\pi_{1F} : F \to E_1$ and $\pi_{2F} : F \to E_2$ to be the restrictions of the projection maps.
Notice that $F$ is a map iff $\pi_{1F}$ is a bijection.  In general, we have:
\begin{equation}\label{2.16}
F \quad = \quad \pi_{2F} \circ (\pi_{1F})^{-1}. \hspace{3cm}
\end{equation}

For $A \times B \subset E_1 \times E_2$ we see that
\begin{equation}\label{2.16a}
\pi_{1F}((A \times B) \cap F) \  = \ A \cap F^{-1}(B) \quad \mbox{and}
\quad \pi_{2F}((A \times B) \cap F) \  = \  F(A) \cap B.
\end{equation}
and so we have
\begin{equation}\label{2.16b}
F(A) \cap B \not= \emptyset \quad \Longleftrightarrow \quad A \cap F^{-1}(B) \not= \emptyset \quad
\Longleftrightarrow \quad (A \times B) \cap F \not= \emptyset.
\end{equation}
Furthermore,
\begin{equation}\label{2.16c}
F(A) \subset B \quad \Longleftrightarrow \quad A \subset F^*(B)\quad \Longleftrightarrow
\quad \pi_{1F}^{-1}(A) \subset A \times B,
\end{equation}
because each of these says
\begin{equation}\label{2.16d}
x \in A \ \mbox{and} \ (x,y) \in F \quad \Longrightarrow \quad y \in B.
\end{equation}

Now we introduce topology.

\begin{prop}\label{prop2.9}(a)  If $f : E_1 \to E_2$ is a continuous map then
the relation $f$ is a closed subset of $E_1 \times E_2$ and $\pi_{1f} : f \to E_1$ is a homeomorphism.

(b) If $f : X_1 \to X_2$ is a map between compact metric spaces which is a closed subset of $X_1 \times X_2$
 then it is a continuous map.
 \end{prop}

 {\bfseries Proof:} (a) If $(x,y) \not\in f$ then there exist disjoint open $U, V \subset E_2$ with
 $f(x) \in U$ and $y \in V$. It follows that $f^{-1}(U) \times V$ is an open subset of $E_1 \times E_2$
 which contains $(x,y)$ and is disjoint from $f$.  As the complement of $f$ is open, $f$ is closed.
 Clearly, $\pi_{1f}$ is continuous and the continuous inverse is given by $x \mapsto (x,f(x))$ for all $x \in E_1$.

 (b) Since $X_1$ and $X_2$ are compact  the closed  $f \subset X_1 \times X_2$ is also a
 compact set.  Hence, the continuous bijection $\pi_{1f} : f \to X_1$ is a homeomorphism. So
 $f = \pi_{2f} \circ (\pi_{1f})^{-1}$ is a continuous map.

 $\Box$ \vspace{.5cm}

For compact metric spaces $X_1, X_2$ a closed (or open) relation $F : X_1 \to X_2$ is a relation which is a closed
(resp. open) subset of $X_1 \times X_2$. For example, if $X$ is a compact metric space and $\epsilon \geq 0$ let
\begin{equation}\label{2.17}
\begin{split}
V_{\epsilon} \quad = \quad \{ (x,y) : d(x,y) < \epsilon \ \} \hspace{2cm} \\
\bar V_{\epsilon} \quad = \quad \{ (x,y) : d(x,y) \leq \epsilon \ \} \hspace{2cm}
\end{split}
\end{equation}
Each $V_{\epsilon}$ is open and each $\bar V_{\epsilon}$ is closed. If $A \subset X$ then $V_{\epsilon}(A)$ is
the $\epsilon$ neighborhood of $A$, i.e. the set of points each of which has distance less than $\epsilon$ from
some point of $A$.

\begin{prop}\label{prop2.10} Let $X_1, X_2, X_3$ be compact metric spaces. Let $F: X_1 \to X_2$ and
$G : X_2 \to X_3$ be relations.
\begin{enumerate}

\item[(a)] If $F$ and $G$ are both closed relations (or both open relations) then the inverse $F^{-1}$
and the composition $G \circ F$ are both closed relations (resp. both open relations).

\item[(b)] Let $A \subset X_1$. If $F$ is an open relation then $F(A)$ is an open subset of $X_2$. If $A$ and
$F$ are both closed then $F(A)$ is closed.

\item[(c)] Let $B \subset X_2$ and let $F$ be a closed relation. If $B$ is closed then $F^{-1}(B)$ is closed.
If $B$ is open then $F^*(B)$ is open.

\item[(d)] Let $F$ be a closed relation and $A$ be a closed subset of $X_1$.  For every $\epsilon > 0$ there
exists $\delta > 0$ such that $ F(V_{\delta}(A)) \subset V_{\epsilon}(F(A))$, equivalently,
$ V_{\delta}(A) \subset F^*(V_{\epsilon}(F(A)))$.

\item[(e)] If $F$ is a closed relation and $A$ is a closed subset of $X_1$ then the restriction
$F \cap (A \times X_2)$ is a closed relation from $A$ to $X_2$.
\end{enumerate}
\end{prop}

{\bfseries Proof:} (a) The results for $F^{-1}$ are obvious since switching coordinates is a homeomorphism
from $X_1 \times X_2$ to $X_2 \times X_1$. $G \otimes F$ is closed (or open) when both $F$ and $G$ are closed
(resp. open).  The projection map $\pi_13$ is always an open map and by compactness it is a closed map. So the results
for $G \circ F$ follow from (\ref{2.12}).

(b) When $F$ is open, each $F(x)$ is open and so $F(A)$ is open by the union representation in (\ref{2.4}).
Notice that no conditions on $A$ are required.  When
$A$ and $F$ are closed then $F(A)$ is closed by the projection representation in (\ref{2.4}).

(c) When $B$ is closed, $F^{-1}(B)$ is closed by (b) applied to $F^{-1}$. Hence, when $B$ is open
$F^*(B) = X_1 \setminus F^{-1}(X_2 \setminus B)$ is open.

(d) Since $F(A) \subset V_{\epsilon}(F(A))$ it follows from (\ref{2.16c}) that the compact set $A$ is contained
in the open set $F^*(V_{\epsilon}(F(A)))$ and hence there exists $\delta > 0$ such that
$V_{\delta}(A) \subset F^*(V_{\epsilon}(F(A)))$.  By (\ref{2.16c}) again it
follows that $ F(V_{\delta}(A)) \subset V_{\epsilon}(F(A))$.

(e) Obvious.

$\Box$ \vspace{.5cm}

If $F: X_1 \to X_2$ and
$G : X_1 \to X_3$ are closed relations we define the product relation $F \bigtriangleup G : X_1  \to X_2 \times X_3$
to be $\{ (x,y,z) : (x,y) \in F, (x,z) \in G  \ \}$.  This is $G \otimes F^{-1}$ with the first and second
 coordinates switched.  Hence, $F \bigtriangleup G$ is a closed relation.

For a closed relation $F : X_1 \to X_2$ we define
\begin{equation}\label{2.18}
ONE_F \quad =_{def} \quad \{ x \in X_1 : F(x) \ \mbox{is a singleton set} \ \}.
\end{equation}
The restriction
\begin{equation}\label{2.19}
f_F \quad =_{def} \quad F \cap (ONE_F \times X_2) \hspace{3cm}
\end{equation}
is a map from $ONE_F$ into $X_2$.

\begin{prop}\label{prop2.11} If $F : X_1 \to X_2$ is a closed relation between compact metric spaces then
$ONE_F$ is a $G_{\delta}$ subset of $X_1$ and $f_F : ONE_F \to X_2$ is a continuous map. \end{prop}

{\bfseries Proof:} For every $\epsilon > 0$ the set $(F \bigtriangleup F)^*(V_{\epsilon})$ consists of the
points $x \in X_1$ such that the diameter of $F(x)$ is less than $\epsilon$.  By Proposition \ref{prop2.10}(c)
this is an open set.  Intersecting over positive rationals we obtain the $G_{\delta}$ set which is the
disjoint union of $ONE_F$ and $F^*(\emptyset)$.  The latter is open and so its complement is closed and hence
$G_{\delta}$. Intersecting with this complement we obtain $ONE_F$ as a $G_{\delta}$.

Now if $x \in ONE_F$, and $\{ (x_n,y_n) \}$ is a sequence in $F$ with $\{x_n \}$ converging to $x$, then for
every limit point $y$ of the sequence $\{ y_n \}$ the point $(x,y) \in F$ because the relation is closed.
By definition of $ONE_F$ it must be that $(x,y) = (x,f_F(x))$.  Since $f_F(x)$ is the only limit point of
the sequence $\{ y_n \}$ and the space $X_2$ is compact it follows that $\{ y_n \}$ converges to $f_F(x)$.
In particular, $f_F$ is continuous on $ONE_F$.

$\Box$ \vspace{.5cm}

Now let $f : X_1 \to X_2$ be a continuous map between compact metric spaces. We say that $f$ is \emph{injective at}
$x \in X_1$ when $f^{-1}(f(x)) = \{x \}$, i.e. $f(x_1) = f(x)$ implies $x_1 = x$. We denote by $IN_f$ the set of
points at which $f$ is injective so that
\begin{equation}\label{2.20}
IN_f \quad =_{def} \quad f^{-1}(ONE_{f^{-1}}) \quad = \quad
\{ x \in X_1 : f^{-1}(f(x)) = \{x \} \ \}.
\end{equation}

For an open $U \subset X_1$ we define the open subset $U^f$ of $U$ by
\begin{equation}\label{2.21}
\begin{split}
U^f \quad = \quad X_1 \setminus f^{-1}(f(X_1 \setminus U)) \quad = \quad f^{-1}(X_2 \setminus f(X_1 \setminus U)) \\
= \quad f^{-1}((f^{-1})^*(U)) \quad = \quad \{ x \in X_1 : f^{-1}(f(x)) \subset U \ \}.
\end{split}
\end{equation}

\begin{prop}\label{prop2.12}If $f : X_1 \to X_2$ is a continuous map between compact metric spaces then
$IN_f$ is a $G_{\delta}$ subset of $X_1$ and $x \in IN_f$ iff for every open $U \subset X_1$, $x \in U$ implies
$x \in U^f$.

The continuous map $f$ restricts to a homeomorphism $f : IN_f \to ONE_{f^{-1}}$.\end{prop}

{\bfseries Proof:}  For a continuous map, $f^{-1}$ is a closed relation and so $ONE_{f^{-1}}$ is a $G_{\delta}$
set by Proposition \ref{prop2.11}. Hence, the pre-image under $f$ is $G_{\delta}$.

Clearly, $f^{-1}(f(x)) = \{ x \}$ implies that $f^{-1}(f(x)) \subset U$ when $x \in U$.  On the other hand,
if $f(x_1) = f(x)$ with $x_1 \not= x$ then $U = X_1 \setminus \{ x_1 \}$ is an open set containing $x$ but not
containing $f^{-1}(f(x))$.

Now it is clear that $f(IN_f) = ONE_{f^{-1}}$. The restriction of $f$ to $IN_f$ is of course continuous. By
Proposition \ref{prop2.11} applied to the closed relation $F = f^{-1}$ the restriction of $f^{-1}$ to $ONE_{f^{-1}}$
defines a continuous map. These are clearly inverse homeomorphisms between the $G_{\delta}$ sets $IN_f \subset X_1$
and $ONE_{f^{-1}} \subset X_2$.

$\Box$ \vspace{.5cm}

Clearly, $f$ is \emph{injective} or \emph{one-to-one} exactly when $IN_f = X_1$.

\begin{df}\label{df2.13} A continuous map $f : X_1 \to X_2$ between compact metric spaces is called
 \emph{almost one-to-one} when $IN_f$ is a dense subset of $X_1$. \end{df}
\vspace{.5cm}

We saw in Proposition \ref{prop2.9}(a) that a closed relation $f : X_1 \to X_2$
between compact metric spaces is a continuous map iff the
projection $\pi_{1f} : f \to X_1$ is a bijection.   For $f$ to be a suitable relation we will weaken this to demand that
$\pi_{1f}$ be a surjection which is almost one-to-one. So we pause to sketch the properties of such maps or, equivalently in this
context, irreducible maps.

\begin{prop} \label{prop2.14} Let $f : X_1 \to X_2$ be a continuous map between compact metric spaces.
The following are equivalent:
\begin{itemize}
\item[(i)]The map $f$ is almost one-to-one.
\item[(ii)] For every nonempty open  $U \subset X_1$ there exists an open $W \subset X_2$ such that
$W \cap f(X_1) \not= \emptyset$ and $f^{-1}(W) \subset U$.
\item[(iii)] $U^f$ is nonempty for every nonempty open
$U \subset X_1$.
\item[(iv)] $U^f$ is dense in $U$ for every open $U \subset X_1$.
\end{itemize}
\end{prop}

{\bfseries Proof:}  (ii) $\Leftrightarrow$ (iii) Clearly, $f^{-1}(W) \subset U$ implies $f^{-1}(W) \subset U^f$.
$f^{-1}(W)$ is nonempty iff $W$ meets $f(X_1)$. Thus, $U^f$ is nonempty if $f^{-1}(W) \subset U$ for a set $W$
which meets $f(X_1)$.  Conversely,  $W = (f^{-1})^*(U)$ is an open subset of $X_2$
with $f^{-1}(W) = U^f \subset U$. If $U^f$ is nonempty then $W$ meets $f(X_1)$.

(i) $\Rightarrow$ (iii) By Proposition \ref{prop2.12} $IN_f \cap U \subset U^f$ for every open set $U$ and so if
$IN_f$ is dense, $U^f$ is nonempty whenever $U$ is.

(iii) $\Rightarrow$ (iv) If $V$ is any open subset of $U$ then $V^f \subset V \cap U^f$.
By (iii) $V^f$ is nonempty whenever $V$ is nonempty and so $U^f$ meets every
nonempty open subset of $U$ and so is dense in $U$.

(iv) $\Rightarrow$ (i) Follow the proof of Proposition \ref{prop2.6}
with $\mathcal{A}$ the same set of pairs as was used there.  We show with a similar argument that
\begin{equation}\label{2.22}
IN_f \quad = \quad \bigcap_{(U_1,U_2) \in \mathcal{A}} \ (U_1)^f \cup (X_1 \setminus \overline{U_2})^f. \hspace{2cm}
\end{equation}
By Proposition \ref{prop2.12} again $IN_f$ is contained in the intersection.
On the other hand if $x$ is in the intersection
and $x \in U$ then choose a pair with $x \in U_2$ and $U_1 \subset U$. Since $x \not\in X_1 \setminus U_2$ it must be
that $x \in U_1^f \subset U^f$.  So by Proposition \ref{prop2.12} $x \in IN_f$.

This expresses $IN_f$ as the intersection of a countable family of dense open sets.  Hence, $IN_f$ is dense by the
Baire Category Theorem.

$\Box$ \vspace{.5cm}

\begin{df}\label{df2.15} Let $f : X_1 \to X_2$ be a  surjective continuous map between compact metric spaces.
A closed subset $A \subset X_1$ is called \emph{minimal for} $f$ when it is minimal in the class of closed subsets of
$X_1$ which are mapped onto $X_2$ by $f$. A continuous map $f : X_1 \to X_2$ is called \emph{irreducible} when it is a
surjection and $X_1$ is minimal for $f$. \end{df}
\vspace{.5cm}

So the continuous surjection $f$ is irreducible when $X_1$ itself is the only closed subset of $X_1$ mapped by $f$ onto $X_2$.  Thus, the restriction of a continuous surjection to a minimal subset is an irreducible map.

\begin{prop}\label{prop2.16} Let $f : X_1 \to X_2$ be a  surjective continuous map between compact metric spaces.
\begin{enumerate}
\item[(a)] $X_1$ contains a closed subset which is minimal for $f$.

\item[(b)]The following are equivalent:
\begin{itemize}
\item[(i)] $f$ is irreducible, i.e. $X_1$ is minimal for $f$.
\item[(ii)] $f$ is almost one-to-one.
\item[(iii)] $D \subset X_1$ is dense in $X_1$ iff $f(D)$ is dense in $X_2$.
\end{itemize}
\end{enumerate}
\end{prop}

{\bfseries Proof:}(a)  By compactness of $X_1$ Zorn's Lemma  applies to the collection of closed subsets
of $X_1$ which are mapped onto $X_2$.  Hence, the collection contains minimal elements.

(b)(ii) $\Rightarrow$ (i)  If $A$ maps onto $X_2$ then clearly $IN_f \subset A$. If $A$ is closed
and $IN_f$ is dense in $X_1$ then $A = X_1$ and so $X_1$ is minimal.

(i) $\Leftrightarrow$ (iii) Since $f$ is surjective, $D$ dense implies $f(D)$ is dense.  On the other hand,
if $f(D)$ is dense then $\overline{D}$ maps onto $\overline{f(D)} = X_2$. If $f$ is irreducible then
$\overline{D}$ must equal $X_1$ and so $D$ is dense. On the other hand, given (iii) $f(A) = X_2$ implies that
$A$ is dense in $X_1$. So if $A$ is closed it equals $X_1$.  Thus, $f$ is irreducible.

(i) $\Rightarrow$ (ii)  Let $U \subset X_1$ be open and nonempty. Since $X_1 \setminus U$ is a proper closed subset
of $X_1$ and $f$ is irreducible, the image $f(X_1 \setminus U)$ is a proper closed subset of $X_2$ and its
complement, $X_2 \setminus f(X_1 \setminus U)$ is a nonempty open subset of $X_2$. Since $f$ is surjective
$U^f = f^{-1}(X_2 \setminus f(X_1 \setminus U))$ is nonempty. By Proposition \ref{prop2.14} $f$ is almost one-to-one.

$\Box$\vspace{.5cm}

\begin{prop}\label{prop2.17} If $f : X_1 \to X_2$ is an irreducible map
then $f$ is almost open with $IN_f = OPEN_f$. \end{prop}

{\bfseries Proof:} If $U$ is a nonempty open subset of $X$ and $x \in U \cap IN_f$
then $x \in U^f = f^{-1}[(f^{-1})^*(U)]$.  Because $f$ is surjective, $f(U^f)$ is the open set
$(f^{-1})^*(U)$ and it contains $f(x)$. Hence, $x \in OPEN_f$.  That is, $IN_f \subset OPEN_f$ and since
$IN_f$ is dense $f$ is almost open.

Now suppose that $f$ is a continuous surjection and that there exists $x \in OPEN_f \setminus IN_f$. We will
show that $f$ is not irreducible. Since $x \not\in IN_f$ there exists a point $x_1$ distinct from $x$ with
$f(x_1) = f(x)$. Let $U, U_1$ be disjoint open sets with $x \in U$ and $x_1 \in U_1$. Since $x \in OPEN_f$,
$f(x) = f(x_1) \in f(U)^{\circ}$. Hence,  $V = U_1 \cap f^{-1}[f(U)^{\circ}]$ is an open set containing
$x_1$ and so $A = X_1 \setminus V$ is a proper closed subset of $X_1$. Note that $f(V) \subset f(U)$ and $U \subset A$.
As $f$ is surjective $X_2 = f(A \cup V) = f(A) \cup f(V) = f(A)$. It follows that $X_1$ is not minimal for $f$.

Thus, if $f$ is irreducible $OPEN_f \subset IN_f$.

$\Box$ \vspace{.5cm}

\begin{prop}\label{prop2.18} Let $f: X_1 \to X_2$ and $g : X_2 \to X_3$ be surjective
continuous maps of compact metric spaces. The composition
$g \circ f : X_1 \to X_3$ is irreducible iff both $f$ and $g$ are irreducible.
\end{prop}

{\bfseries Proof:} The proof is an easy exercise using minimality.

%

$\Box$ \vspace{.5cm}

\begin{prop}\label{prop2.19} Let $f : X_1 \to X_2$ be a  surjective continuous map between compact metric spaces.
The following are equivalent:
\begin{itemize}
\item[(i)] $X_1$ contains a unique closed subset which is minimal for $f$.
\item[(ii)] $f(IN_f) = ONE_{f^{-1}}$ is dense in $X_2$.
\end{itemize}
When these conditions hold the unique minimal set is the closure $\overline{IN_f}$.
\end{prop}

{\bfseries Proof:} If $A \subset X_1$ maps onto $X_2$ then $A$ contains $IN_f$ and so if $A$ is closed
it contains $\overline{IN_f}$.  On the other hand, $f(IN_f)$ is dense then $\overline{IN_f}$ maps onto
$X_2$.  Thus, condition (ii) implies that $\overline{IN_f}$ is the unique closed subset of $X_1$ which is minimal
for $f$.

To complete the proof we need to use some results from Akin (1993)

\begin{lem}\label{lem2.20} If $f : X_1 \to X_2$ is a surjective continuous map between compact metric spaces
then $(f^{-1})^*(OPEN_f) = \{ y \in X_2 : f^{-1}(y) \subset OPEN_f \}$ is a dense $G_{\delta}$ subset of
$X_2$. \end{lem}

{\bfseries Proof:} The closed relation $f^{-1}: X_2 \to X_1$ is lower semicontinuous at $y \in X_2$ when
$U \subset X_1$ is open and $f^{-1}(y) \cap U \not= \emptyset$ implies $\{ y_1 : f^{-1}(y_1) \cap U \not= \emptyset \}$
is a neighborhood of $y$.  This says that for all $x \in f^{-1}(y)$ and $U$ open with $x \in U$, $f(U)$ is a
neighborhood of $y = f(x)$.  That is, $f$ is open at every point of $f^{-1}(y)$.  See Akin (1993) Proposition 7.11.
By Theorem 7.19 of Akin (1993)  the set of points $y \in X_2$ at which $f^{-1}$ is lower semicontinuous is a dense
$G_{\delta}$ subset of $X_2$. This is $(f^{-1})^*(OPEN_f)$.

$\Box$ \vspace{.5cm}

Now assume that $f(IN_f)$ is not dense in $X_2$.  By Lemma \ref{lem2.20} we can choose a point $y$ in the complement
of the closure of $f(IN_f)$ such that $f^{-1}(y) \subset OPEN_f$. Since $y \not\in f(IN_f)$ we can choose  open
 $U_1, U_2 \subset X_1$  with disjoint closures such that each meets $f^{-1}(y)$.
 Let $V = f(U_1)^{\circ} \cap f(U_2)^{\circ}$ which is an open set containing $y$ because $f$ is open at
 the points of $f^{-1}(y)$.  For $i = 1,2$ let $A_i = (X_1 \setminus f^{-1}(V)) \cup \overline{U_i}$. Each of these is
 a proper closed subset of $X_1$.  For example, $A_1$ is disjoint from $U_2 \cap f^{-1}(V)$. On the other hand,
 since $f$ is surjective and $V \subset f(U_1) \cap f(U_2)$ we have $X_2 = f(A_i \cup f^{-1}(V)) = f(A_i) \cup V =
 f(A_i)$.  By Zorn's Lemma each $A_i$ contains a minimal set $M_i$ for $f$.  They cannot be the same set
 because $M_1 \cap M_2 \subset A_1 \cap A_2 = X_1 \setminus f^{-1}(V)$
 whose image under $f$ does not contain the points
 of $V$. Thus, there are at least two distinct subsets minimal for $f$.

 This completes the proof that (i) implies (ii).

 $\Box$ \vspace{.5cm}

 Recall that a point $x$ is called \emph{isolated} when it is open as well as closed.

 \begin{lem}\label{lem2.21} Let $f : E_1 \to E_2$ be a continuous map and let $x \in E_1$
 \begin{enumerate}
 \item[(a)] If $x$ is isolated in $E_1$ and $f$ is weakly open then $f(x)$ is is isolated in $E_2$.

 \item[(b)] If $f(x)$ is isolated in $E_2$ and $IN_f$ is dense in $E_1$ then $x$ is isolated in $E_1$ and $x \in IN_f$.

 \item[(c)] If $f$ is an irreducible map between compact metric spaces then $x$ is isolated iff $f(x)$ is
 isolated and in that case, $x \in IN_f$.
 \end{enumerate}
 \end{lem}

 {\bfseries Proof:} (a) If $x$ is open and $f$ is weakly open then the interior of $\overline{f(x)} = f(x)$ is nonempty which implies that
 $f(x)$ is open.

 (b) If $f(x)$ is open and $IN_f$ is dense then the nonempty open set $f^{-1}(f(x))$ meets $IN_f$. If $x_1$ is in the intersection then
 $\{ x_1 \} = f^{-1}(f(x_1)) = f^{-1}(f(x))$. As $x$ is in this set, we have $x_1 = x$. Since $f(x)$ is open, $\{ x \} = f^{-1}(f(x))$ is open.

 (c) An irreducible map is both almost open and almost one-to-one by Propositions \ref{prop2.16} and \ref{prop2.17}.


 $\Box$ \vspace{.5cm}

 Recall the Uniqueness of Cantor Theorem which says that any compact,
 zero-dimensional metric space without isolated points is homeomorphic to
 the Cantor Set.

\begin{theo}\label{theo2.22} If $X$ is a compact metric space without
isolated points then there exists an irreducible map $f : C \to X$
 with $C$ the Cantor Set. If, in addition, $X$ is connected then $IN_f$ has
 empty interior.  That is, $C \setminus IN_f$ is dense in $C$. \end{theo}

 {\bfseries Proof:} It is well known that there exists a continuous map $F$ from a Cantor set $C$
 onto any compact metric space $X$.

%

By Proposition \ref{prop2.16}(a) there is a closed subset $Z $ of $C$ which is
minimal for $F$. Since $Z$ is a closed subset of $C$ it is
a compact, zero-dimensional metric space.  By Lemma \ref{lem2.21} it has no isolated points
since $X$ does not and since $F : Z \to X$ is
irreducible.  Hence, there is a homeomorphism $h : C \to Z$. Let $ f = F \circ h$.

If $(IN_f)^{\circ}$ is nonempty then because $C$ is zero-dimensional there is a  clopen  $A \subset IN_f$ which
is neither empty nor all of $C$. Since $f$ is irreducible $f(A)$ is a nonempty, proper, closed subset of $X$.
By  Proposition \ref{prop2.17} $f$ is open at every point of $A \subset IN_f$. Hence, $f(A)$ is open and so
$X$ is not connected.

 $\Box$ \vspace{1cm}

 From this we notice a problem that might arise with  selection functions.

 \begin{theo}\label{theo2.23} Let $f : X_2 \to X_1$ be an irreducible map. If $g : X_1 \to X_2$
 is a selection function for the relation $f^{-1} : X_1 \to X_2$, i.e. $g$ is a map and $g \subset f^{-1}$,
then $(IN_f)^{\circ}$ is a dense subset of $g(X_1)^{\circ}$. In particular, if $(IN_f)^{\circ} = \emptyset$
then $g(X_1)^{\circ} = \emptyset$.
\end{theo}

{\bfseries Proof:} If $f^{-1}(f(x)) = \{ x \}$ then $g(f(x)) = x$. Hence, $IN_f \subset g(X_1)$ and so
we have inclusion of the interiors as well.

Now let $U$ be any nonempty open subset disjoint from $(IN_f)^{\circ}$. We will show that $U \setminus g(X_1)$
is nonempty. This will show that $g(X_1)^{\circ}$ is contained in the closure of $(IN_f)^{\circ}$.

Because $f$ is irreducible, it is almost one-to-one by Proposition \ref{prop2.16} and so by Proposition
\ref{prop2.14} there is a nonempty open set $V \subset X_1$ such that $f^{-1}(V) \subset U$.
Since $U$ is disjoint from  $(IN_f)^{\circ}$ the open set $f^{-1}(V)$ meets the complement of $IN_f$.
That is, there exist distinct points $x_1, x_2$ with $x_1 \in f^{-1}(V)$ and with a common value $y$.
That is, $y = f(x_1) = f(x_2)$. Because
$x_1 \in f^{-1}(V), \ y \in V$ and so $x_2 \in f^{-1}(V)$ as well. So either $x_1 \not= g(y)$ or
$x_2 \not= g(y)$. Hence, one or the other is a point of $U \setminus g(X_2)$.

$\Box$ \vspace{.5cm}

{\bfseries Remark:}  $(IN_f)^{\circ}$ might be a proper subset of $g(X_1)^{\circ}$. For example, let
 $X_2 = [ - \frac{1}{2},0] \cup [1,\frac{3}{2}]$, $X_1 = [0,1]$ and define $f :X_2 \to X_1$ by
\begin{equation}\label{3.9a}
f(x) = \begin{cases} \ x + \frac{1}{2} \qquad \mbox{for} \quad - \frac{1}{2} \leq x \leq 0, \\
 \ x - \frac{1}{2} \qquad \mbox{for} \quad 1 \leq x \leq \frac{3}{2}. \end{cases}
 \end{equation}
 Clearly, $f$ is irreducible with $IN_f$ the open set $[\frac{1}{2},0) \cup (1,\frac{3}{2}]$
 and if $g$ is either of the two selection functions for $f^{-1}$ then $g(X_1)$ is an open set containing
 the endpoint $g(\frac{1}{2})$ as well.
 \vspace{.5cm}

Thus, if $X$ is a compact connected metric space with more than one point, and so without isolated points, then
there exists an irreducible map $f : C \to X$ and if $g : X \to C$ is any selection function for $f^{-1}$ then
$g(X)$ has empty interior in $C$.  By  Theorem \ref{theo1.1aa} the map $g$ is quasi-continuous but it is as far
as possible from being an almost open map.

Furthermore, if $f_1 : C \to X$ and $f_2 : C \to X$ are irreducible maps then $F = f_2 \circ (f_1)^{-1} : X \to X$
is a closed relation with $\pi_{1F}$ irreducible. In fact $F$ is a suitable relation.
  If $g_1 : X \to C$ is a selection function
for $(f_1)^{-1}$ then $f_2 \circ g_1 : X \to X$ is a selection function for $F$ and so is quasi-continuous.
If $IN_{f_1} \cup IN_{f_2} = C$,
or, equivalently the complements of $IN_{f_1}$ and $IN_{f_2}$ are disjoint, then $(f_2 \circ g_1)(X))^{\circ}
= \emptyset$.  To prove this let $U \subset X$ be an arbitrary nonempty  open set.
Then $ (f_2)^{-1}(U)$ is a nonempty open subset of $C$ and so contains
an $(f_1)^{-1}(V)$ for some nonempty open $V \subset X$.
As in the above proof there exist distinct $x_1, x_2 \in C$ with
$y = f_1(x_1) = f_1(x_2) \in V$. Suppose $g_1(y) \not= x_1$.  Since $x_1 \not\in IN_{f_1}$ it follows that
$x_1 \in IN_{f_2}$ and so there does not exist any $x \not= x_1$ with $f_2(x) = f_2(x_1)$.  In particular,
there does not equal any $z \in X$ such that $f_2(g(z)) = f_2(x_1)$.  Hence, $f_2(x_1) \in U \setminus f_2(g_1(X))$.
\vspace{.5cm}

{\bfseries Example:}  Let $f_1 : C = \{ 0,1 \}^{\mathbb N} \to [0,1]$ be the map given by
$z \mapsto \Sigma_i \frac{z_i}{2^{i+1}}$
then the complement of $IN_{f_1}$ is contained in the countable set of sequences
$z$ such either $z_i $ eventually $0$ or is
eventually $1$. Let $\tilde{0} = 1 $ and $\tilde{1} = 0$. Define the homeomorphism
$h $ on $C$ by
\begin{equation}\label{3.10a}
h(z)_i \quad = \quad \begin{cases} \quad z_i \quad \mbox{for} \ i \ \mbox{even} \\
\quad \tilde{z_i} \quad \mbox{for} \ i \ \mbox{odd}. \end{cases}
\end{equation}
Let $f_2 = f_1 \circ h$. The complement of $IN_{f_2}$ is contained in the countable set of sequences which are
eventually $010101...$. This is disjoint from the complement of $IN_{f_1}$.
\vspace{1cm}

 \section{Suitable Relations}

Let $F : X_1 \to X_2$ be a closed relation between compact metric spaces.  For $i = 1,2$ we denote
by $\pi_{iF} : F \to X_i$ the restrictions of the projections.  We have
\begin{equation}\label{3.1}
F \quad = \quad \pi_{2F} \circ (\pi_{1F})^{-1}.
\end{equation}

\begin{lem}\label{lem3.1} The following are equivalent:
\begin{itemize}
\item[(i)] $F$ has full domain, i.e. $F^{-1}(X_2) = X_1$.
\item[(ii)] For every $x \in X_1$, $F(x)$ is nonempty.
\item[(iii)] $\pi_{1F}$ is surjective.
\item[(iv)] $F^*(\emptyset) = \emptyset$.
\item[(v)]  For some $B \subset X_2$, $F^*(B) \subset F^{-1}(B)$.
\item[(vi)]  For all $B \subset X_2$, $F^*(B) \subset F^{-1}(B)$.
\end{itemize}
\end{lem}

{\bfseries Proof:} An easy exercise applying the various definitions.

$\Box$ \vspace{.5cm}

For $(x,y) \in F$ it is clear that $\pi_{1F}$ is injective at $(x,y)$, i.e. $(x,y) \in IN_{\pi_{1F}}$,
precisely when
$F(x) = \{ y \}$ and so when $x \in ONE_F$.

\begin{prop}\label{prop3.2} $\pi_{1F}$ restricts to a homeomorphism from $IN_{\pi_{1F}} \subset F$ onto
$ONE_F \subset X_1$. With $f_F : ONE_F \to X_2$ the continuous map obtained by restricting the closed relation
$F$ to $ONE_F$ we see that
\begin{equation}\label{3.2}
f_F \quad = \quad IN_{\pi_{1F}} \quad \subset \quad ONE_F \times X_2. \hspace{2cm}
\end{equation}
\end{prop}

{\bfseries Proof:} That $\pi_{1F}$ restricts to a homeomorphism follows from Proposition \ref{prop2.12}
applied to the map $\pi_{1F}$.  Notice that $ONE_F = ONE_{(\pi_1F)^{-1}}$.
The set $IN_{\pi_{1F}}$ is just the restriction of $F$ to $ONE_F$ and
this is the  map $f_F$ which is continuous by Proposition \ref{prop2.11}.

$\Box$ \vspace{.5cm}

 \begin{cor}\label{cor3.2a} Let $F : X_1 \to X_2$ be a closed relation between compact metric spaces.
The following are equivalent:
\begin{itemize}
\item[(i)] The map $\pi_{1F} : F \to X_1$ is surjective and $F$ contains a unique
closed subset which is minimal for $\pi_{1F}$.
\item[(ii)] $ONE_F$ is dense in $X_1$.
\end{itemize}
When these conditions hold the unique minimal set is the closure in $X_1 \times X_2$ of the map $f_F : ONE_F \to X_2$.
\end{cor}

{\bfseries Proof:} The closed set $F^{-1}(X_2)$ contains the set $ONE_F$ and so it equals $X_1$
when the latter is dense. So by Lemma \ref{lem3.1}
$\pi_{1F}$ is surjective when $ONE_F$ is dense.  The result then follows from
Proposition \ref{prop2.19} and Proposition \ref{prop3.2}.

$\Box$ \vspace{.5cm}

By  Proposition \ref{prop2.16}(a) any closed full domain relation contains at least one minimal closed full domain relation.

\begin{theo}\label{theo3.3} Let $F: X_1 \to X_2$ be a closed relation between compact metric spaces
and $\pi_{1F} : F \to X_1$ be the restriction to $F \subset X_1 \times X_2$ of the first coordinate projection map.
\begin{enumerate}
\item[(a)] The following conditions are equivalent.
\begin{itemize}
\item[(i)] $F$ is a minimal closed full domain relation from $X_1$ to $X_2$.
\item[(ii)] The map $\pi_{1F}$ is irreducible, i.e. $F$ is minimal for $\pi_{1F}$.
\item[(iii)] The map $\pi_{1F}$ is surjective and almost one-to-one.
\item[(iv)] $F$ has full domain and for every open $U \times V \subset X_1 \times X_2$ which meets $F$ there exists
$x \in U$ such that $F(x) \subset V$.  That is, $U \cap F^*(V)$ is nonempty.
\item[(v)] For every open $V \subset X_2 \ F^*(V)$ is a dense subset of $F^{-1}(V)$.  That is,
\begin{equation}\label{3.3}
F^*(V) \quad \subset \quad F^{-1}(V) \quad \subset \quad \overline{F^*(V)}.
\end{equation}
\end{itemize}

\item[(b)]  If  $F$ is a minimal full domain closed relation then $ONE_F$ is dense in $X_1$ and $F$ is the closure in $X_1 \times X_2$
of the map $f_F : ONE_F \to X_2$. That is, $f_F$ is dense in $F$.

\item[(c)] If there exists  $D $ a dense subset of $ X_1$ and a continuous map $f : D \to X_2$ such that
$f$ is dense in $F$ then $\pi_{1F}$ is irreducible, $D \subset ONE_F$ and $f = f_F$ on $D$.

\end{enumerate}
\end{theo}

{\bfseries Proof:} (a)(i) $\Leftrightarrow$ (ii) This is the definition of an irreducible map.

(ii) $\Leftrightarrow$ (iii) Apply Proposition \ref{2.16}.

(iii) $\Rightarrow$ (iv) $X_1 = F^{-1}(X_2)$ because $\pi_{1F}$ is surjective. If an open $U \times V$ meets
$F$ then $F \cap (U \times V)$ is a nonempty open subset of $F$. Because $\pi_{1F}$ is almost one-to-one
Proposition \ref{prop2.14} implies there exists $x \in X_1$ such that $(\pi_{1F})^{-1}(x) \subset U \times V$
or, equivalently, $x \in U$ and $F(x) \subset V$.

(iv) $\Rightarrow$ (v) Because $F$ has full domain Lemma \ref{lem3.1} implies $F^*(V) \subset F^{-1}(V)$.
If $x \in F^{-1}(V)$ then there exists $y \in V$ with $(x,y) \in F$. For any open set $U$ with $x \in U$,
(iv) implies there exists $x_1 \in U$ such that $F(x_1) \subset V$. Thus, $U$ meets $F^*(V)$ and so the latter is
dense in $F^{-1}(V)$.

(v) $\Rightarrow$ (iii) By Lemma \ref{lem3.1} $\pi_{1F}$ is surjective because $F^*(V) \subset F^{-1}(V)$.
Any open subset of $F$ contains a nonempty open set of the form $F \cap (U \times V)$, i.e. $U$ meets $F^{-1}(V)$.
Since $U$ is open (v) implies that $U$ meets $F^*(V)$. There is a point $x \in U$ such that $F(x) \subset V$
and so $\pi_{1F}^{-1}(x) \subset F \cap (U \times V)$.  It follows from Proposition \ref{prop2.14} that $\pi_{1F}$
is almost one-to-one.

(b) Since $\pi_{1F}$ is almost one-to-one, $IN_{\pi_{1F}}$ is dense in $F$. By (\ref{3.2}) this says $f_F$ is dense
in $F$. Since $\pi_{1F}$ is surjective, it follows that $ONE_F = \pi_{1F}(IN_{\pi_{1F}})$ is dense in $X_1$.

(c) If $(x,y) \in F$ with $x \in D$ then since $f$ is dense in $F$ there is a sequence $\{ x_n \}$ in $D$
such that $\{ (x_n,f(x_n)) \}$ converges to $(x,y)$. By continuity of $f$ on $D$ $\{ f(x_n) \}$ converges to
$f(x)$ and so $y = f(x)$.  Thus, for $x \in D$, $F(x) = \{ f(x) \}$. It follows that $D \subset ONE_F$ and that on
$D$ the map $f_F$ restricts to $f$. Finally, $IN_{\pi_{1F}} = f_F$ contains $f$ and so is dense in $F$.  Hence,
$\pi_{1F}$ is almost one-to-one.  Since $D$ is dense in $X_1$, $\pi_{1F}(F)$ contains the closure of $D$ which is
$X_1$.  That is, $\pi_{1F}$ is surjective.

$\Box$ \vspace{.5cm}

{\bfseries Remarks:}  In part (a) it does not suffice that $F^*(V)$ be nonempty for all open $V \subset X_2$
which meet $F(X_1)$.  Let $X_1$ be the disjoint union of two nontrivial spaces $X$ and $Y$ and let $X_2 = X$.
For $F = 1_X \cup (Y \times X)$ it follows that for every proper open subset $V \subset X_2 = X$ that
$F^*(V) = V \subset X \subset X_1$ but for no point $x \in Y$ is $F(x) \subset V$.

\vspace{.5cm}

\begin{cor}\label{cor3.3a} Let $F : X_1 \to X_2$ be  a closed relation between compact metric spaces and
let $f : D \to X_2$ be a continuous map in with $D \subset X_1$. The following are equivalent:
\begin{itemize}
\item[(i)] $D$ is dense in $X_1$ and $F$ is the closure in $X_1 \times X_2$ of $f$.

\item[(ii)] $F$ has full domain and $F$ is the closure in $X_1 \times X_2$ of $f$.

\item[(iii)] $F$ is a minimal closed full domain relation, $D$ is a
dense subset of $ONE_F$ and $f$ is the restriction of
$f_F$ to $D$.
\end{itemize}
When these conditions hold we will say \emph{ the map $f$ is dense in $F$ }.
\end{cor}

{\bfseries Proof:}  (i) $\Leftrightarrow$ (ii) If $f$ is dense in $F$ then $D = \pi_{1F}(f) $ is dense
in $\pi_{1F}(F)$. Then $D$ is dense in $X_1$ iff  $F^{-1}(X_2) = \pi_{1F}(F) = X_1$.

(i) $\Rightarrow$ (iii) Apply Theorem \ref{theo3.3}(c).

(iii) $\Rightarrow$ (i) Assuming (iii)
Theorem \ref{theo3.3}(b) says that $f_F : ONE_F \to X_2$ is dense in $F$.  If  $x \in ONE_F$ then it is
the limit of some sequence $\{ x_n \}$ in $D$ because $D$ is dense on $ONE_F$. By continuity
of $f_F$, $f_F(x)$ is the limit of the sequence $\{ f(x_n) = f_F(x_n) \}$. Thus, $f$ is dense in $f_F$ and so
is dense in $F$.  $D$ is dense in $X_1$ because it is dense in $ONE_F$ and the latter is dense in $X_1$.

$\Box$ \vspace{.5cm}

While proper containment between closed relations can occur, we have

\begin{lem}\label{lem3.3b} Let $\tilde F, F : X_1 \to X_2$ be closed relations such that $\tilde F$ has full domain and
$\pi_{1F}$ is irreducible.  If $\tilde F \subset F$ then $\tilde F = F$. \end{lem}

{\bfseries Proof:}  $\tilde F$ is a closed full domain relation contained in $F$ and $F$ is a minimal closed
full domain relation.  Hence, $\tilde F = F$.
%

$\Box$ \vspace{.5cm}

\begin{theo}\label{theo3.4} Let $F : X_1 \to X_2$ be a closed relation between compact metric spaces with
$\pi_{1F} : F \to X_1$ irreducible. The following conditions are equivalent.
\begin{itemize}
\item[(i)]  $\pi_{2F} : F \to X_2$ is an almost open continuous map.

\item[(ii)] $f_F : ONE_F \to X_2$ is a weakly open continuous map.

\item[(iii)] If $U \subset X_1$ is open and nonempty then $F(U)^{\circ} \not= \emptyset$.

\item[(iv)] If $U \subset X_1$ is open then the open subset
\begin{equation}\label{3.4}
 U_F \quad =_{def} \quad U \cap F^*( F(U)^{\circ} \ )
 \end{equation}
 is dense in $U$.

 \item[(v)] If $D$ is a dense subset of $X_2$ then $F^{-1}(D)$ is dense in $X_1$.

 \item[(vi)] If $D$ is a dense open subset of $X_2$ then $F^*(D)$ is a dense open subset of $X_1$.

 \item[(vii)] If $D$ is a dense $G_{\delta}$ subset of $X_2$ then $F^*(D)$ is a dense $G_{\delta}$ subset of $X_1$.

\item[(viii)] If $B$ is a closed nowhere dense subset of $X_2$ then $F^{-1}(B)$ is a closed nowhere dense subset of
$X_1$.
\end{itemize}
\end{theo}

{\bfseries Proof:} (i) $\Leftrightarrow$ (ii) By Theorem \ref{theo2.3} $\pi_{2F}$ is almost open iff it is
weakly open and by the Variation of Domain  Theorem \ref{theo2.4} it is weakly open iff its restriction to a dense set is weakly open.
Since $\pi_{1F}$ is irreducible, $f_F$ is a dense subset of $F$.  Furthermore, $\pi_{1F}$ restricts to a homeomorphism
of $f_F$ onto its domain $ONE_F$. Thus, the restriction of $\pi_{2F}$ to $f_F$ is weakly open iff
$f_F = (\pi_{2F}|f_F) \circ (\pi_{1F})^{-1}|ONE_F$ is weakly open.

(i) $\Rightarrow$ (iii) If $U \subset X_1$ is open and nonempty then the open set $U \times X_2$ meets $F$ because
$\pi_{1F}$ is surjective. Since $\pi_{2F}$ is almost open $F(U) = \pi_{2F}(F \cap (U \times X_2))$ has a nonempty
interior.

(iii) $\Rightarrow$ (iv)  It suffices as usual to show that $U_F$ is nonempty when $U$ is because
we then apply the result to arbitrary open subsets $V$ of $U$ and note that $V_F \subset V \cap U_F$.
By (iii) $F(U)^{\circ}$ is nonempty and by (\ref{3.3}) $F^*(F(U)^{\circ})$ is dense in $F^{-1}(F(U)^{\circ})$.
Since $U$ is open it suffices to show that $U$ meets $F^{-1}(F(U)^{\circ})$.  But if $y \in F(U)^{\circ}$
then since $y \in F(U)$ there exists $x \in U$ such that $y \in F(x)$.  Hence, $x \in F^{-1}(F(U)^{\circ})$.

(iv) $\Rightarrow$ (iii) Obvious because $F^*(\emptyset) = \emptyset$ by Lemma \ref{lem3.1}.

(iii) $\Rightarrow$ (v) If $U \subset X_1$ is open an nonempty then $F(U)$ has a nonempty interior by (iii) and
so it meets the dense set $D$.  Hence, $U$ meets $F^{-1}(D)$ by (\ref{2.16b}). As $U$ was arbitrary $F^{-1}(D)$ is
dense.

(v) $\Rightarrow$ (vi) $D$ open implies $F^*(D)$ is open and by (\ref{3.3}) it is dense in $F^{-1}(D)$ which is
dense in $X_1$ by (v).

(vi) $\Leftrightarrow$ (vii) If $D$ is open the $F^*(D)$ is open and so by (\ref{2.10}) if $D$ is $G_{\delta}$ then
$F^*(D)$ is $G_{\delta}$. So (vii) obviously implies (vi) and (vi) implies (vii) by the Baire Category Theorem.

(vi) $\Leftrightarrow$ (viii) If $D$ is the complement of $B$ then $F^*(D)$ is the complement of $F^{-1}(B)$.
An open set is dense iff its complement is nowhere dense.

(vi) $\Rightarrow$ (ii) Let $D$ be an open dense subset of $X_2$. By (vi), $F^*(D)$ is an open dense subset of
$X_1$.  Clearly, $ONE_F \cap F^*(D) = (f_F)^{-1}(D)$. Since $ONE_F$ is dense in $X_1$ it follows from Lemma \ref{lem2.5}
(a) that $ONE_F \cap F^*(D)$ is dense in $F^*(D)$ and hence in $X_1$. So as a subset of $ONE_F$ it is dense in $ONE_F$.
Thus, $D$ open and dense in $X_2$ implies $(f_F)^{-1}(D)$ is open and dense in $ONE_F$.  It follows from Proposition
\ref{prop2.2a} that $f_F$ is weakly open.

$\Box$ \vspace{.5cm}

Notice that if $F \subset [0,1] \times [0,1]$ is the closed full domain relation  defined by:
\begin{equation}\label{3.4aa}
F \quad = \quad [0,1] \times \{ 0 \} \cup \{ (\frac{m}{n},y) : y \leq \frac{1}{n} \quad \mbox{for all} \quad m \leq n  \},
\end{equation}
where $m$ and $n$ vary over the positive integers, then $\pi_{2F}$ is almost open. However, the unique minimal full domain
relation contained in $F$ is the constant function $[0,1] \times \{0\}$ which is not almost open.

\begin{df}\label{df3.5} A \emph{suitable relation} is a closed relation $F : X_1 \to X_2$ between compact metric
spaces such that $\pi_{1F} : F \to X_1$ is irreducible and $\pi_{2F} : F \to X_2$ is almost open. \end{df}
\vspace{.5cm}

\begin{prop}\label{prop3.6} (a)Let $f : X_1 \to X_2$ be a continuous map between compact spaces. $f$ is a suitable
relation iff it is an almost open map.

(b) Assume $X_1$ and $X_2$ are compact metric spaces and that
 $f : D \to X_2$ is a continuous map with $D$ dense in $X_1$. Let $F$ be the closure of $f$ in $X_1 \times X_2$.
 $F$ is a suitable relation iff $f$ is a weakly open map. \end{prop}

{\bfseries Proof:} (a) $\pi_{1f} : f \to X_1$ is a homeomorphism and so is irreducible.
Since $f = \pi_{2f} \circ (\pi_1f)^{-1}$ the map $f$ is almost open iff the map $\pi_{2f}$ is almost open.

(b) By Theorem \ref{theo3.3}(c) $\pi_{1F}$ is irreducible and $f$ is the restriction to $D \subset ONE_F$ of
$f_F$. $\pi_{2F}$ is almost open iff $f_F$ is weakly open by Theorem \ref{theo3.4} and so iff $f$ is weakly open by
the Variation of Domain  Theorem \ref{theo2.4}.

$\Box$ \vspace{.5cm}

%
%
%
%

Now we consider composition of suitable relations.

\begin{lem}\label{lem3.8} If $F : X_1 \to X_2$ and $G : X_2 \to X_3$ are closed relations between compact
metric spaces such that  $X_2 = G^{-1}(X_3)$ and $X_1 = F^{-1}(X_2)$, then
\begin{equation}\label{3.5}
ONE_{G \circ F} \quad \subset \quad F^*(ONE_G) \quad \subset \quad F^{-1}(ONE_G), \hspace{2cm}
\end{equation}
and
\begin{equation}\label{3.6}
\begin{split}
ONE_F \cap ONE_{G \circ F} \quad = \quad ONE_F \cap F^*(ONE_G)  \hspace{3cm}\\
= \quad ONE_F \cap F^{-1}(ONE_G) \quad = \quad (f_F)^{-1}(ONE_G). \hspace{1cm}
\end{split}
\end{equation}
\end{lem}

{\bfseries Proof:}  If $y \in F(x)$ then $G(y) \subset (G \circ F)(x)$ with equality if $F(x) = \{ y \}$.
In particular, if $F(x)$ contains some point $y$ not in $ONE_G$ then $x$ is not in $ONE_{G \circ F}$ and
conversely if $x \in ONE_F$.  Thus, $ONE_{G \circ F} \subset F^*(ONE_G)$ with equality after intersection by
$ONE_F$.  Since $X_1 = F^{-1}(X_2)$ the second inclusion in (\ref{3.5}) follows from Lemma \ref{lem3.1}.
For any  $B \subset X_2$ we clearly have
\begin{equation}\label{3.7}
ONE_F \cap F^*(B) \quad = \quad ONE_F \cap F^{-1}(B) \quad = \quad (f_F)^{-1}(B).
\end{equation}
From this follow the remaining equalities in (\ref{3.6}).

$\Box$\vspace{.5cm}

We can restrict the map $f_F : ONE_F \to X_2$ to $(f_F)^{-1}(ONE_G) \subset ONE_F$ and obtain a map
from $(f_F)^{-1}(ONE_G)$ to $ONE_G$ which can be composed with $f_G :ONE_G \to X_3$.  That is, we have
the composed map  $f_G \circ f_F : (f_F)^{-1}(ONE_G) \to X_3$.

\begin{theo}\label{theo3.9} Assume that $F : X_1 \to X_2$ and $G : X_2 \to X_3$ are suitable relations.
The set $(f_F)^{-1}(ONE_G)$ is a $G_{\delta}$ subset of $X_1$ which is dense in $X_1$. It is contained in
$ONE_{G \circ F}$ and on it the functions $f_G \circ f_F$ and $f_{G \circ F}$ agree. The function
$f_{G \circ F}: ONE_{G \circ F} \to X_3$ and its restriction $f_G \circ f_F : (f_F)^{-1}(ONE_G) \to X_3$ are both
weakly open maps.\end{theo}

{\bfseries Proof:} Because $G$ is suitable, $ONE_G$ is a dense $G_{\delta}$ subset of $X_2$. Since $F$ is suitable
$f_F$ is a weakly open map and so by Proposition \ref{prop2.2a} $(f_F)^{-1}(ONE_G)$ is a dense $G_{\delta}$
subset of $ONE_F$. Because $F$ is suitable the latter is a dense $G_{\delta}$ in $X_1$ it follows that $(f_F)^{-1}(ONE_G)$ is
a dense $G_{\delta}$ in $X_1$. By (\ref{3.6}) $(f_F)^{-1}(ONE_G) \subset ONE_{G \circ F}$ and for $x \in (f_F)^{-1}(ONE_G)$
it is clear that the single point in $G \circ F(x)$ is the point$ f_G \circ f_F(x)$.  That is, $f_{G \circ F}$
agrees with $ f_G \circ f_F $ on $(f_F)^{-1}(ONE_G)$.

Because $F$ and $G$ are suitable the maps $f_F : ONE_F \to X_2$ and $f_G : ONE_G \to X_3$ are weakly open. By
the Variation of Domain  Theorem \ref{theo2.4} the restriction  $f_F : (f_F)^{-1}(ONE_G) \to ONE_G$ is weakly open.  By Corollary
\ref{cor2.2b} the composition $f_G \circ f_F : (f_F)^{-1}(ONE_G) \to X_3$ is weakly open and by
the Variation of Domain  Theorem \ref{theo2.4} again the extension to  $f_{G \circ F} : ONE_{G \circ F} \to X_3$ is weakly open as well.

$\Box$ \vspace{.5cm}

\begin{lem}\label{lem3.10} If $F : X_1 \to X_2$ and $G : X_2 \to X_3$ are suitable relations then
$f_G \circ f_F \subset (f_F)^{-1}(ONE_G) \times X_3$ is a dense subset of
\begin{displaymath}
\pi_{13}((X_1 \times ONE_G \times X_3) \cap G \otimes F)
\end{displaymath}
where $\pi_{13} : X_1 \times X_2 \times X_3 \to
X_1 \times X_3$ is the projection to the first and third coordinates. \end{lem}

{\bfseries Proof:}  Recall that $G \otimes F = \{ (x,y,z) : (x,y) \in F$ and $(y,z) \in G \}$.
For any $(x,y,z) \in G \otimes F$ we can choose a sequence $\{ (x_n,y_n) \in (f_F) \cap (f_F)^{-1}(ONE_G) \times X_2$
converging to $(x,y)$ because the map $f_F$ is dense in $F$ and  $(f_F)^{-1}(ONE_G)$ is dense in
$ONE_F$. See Corollary \ref{cor3.3a}.  Since $y_n$ must equal $f_F(x_n)$  and $x_n \in (f_F)^{-1}(ONE_G)$
it follows that the sequence $\{ y_n \}$ is in $ONE_G$. Let $z_n  \in X_3 $ be the unique point in $G(y_n)$.
By going to a subsequence if necessary we can assume that $\{ z_n \}$ converges to a point $\tilde z \in X_3$.
Since $G \otimes F$ is closed we have that $(x,y,\tilde z) \in G \otimes F$.  Now if $y \in ONE_G$ then
$z, \tilde z \in G(y)$ implies $z = \tilde z$. Thus, when $(x,y,z) \in G \otimes F$
with $y \in ONE_G$ we obtain a sequence $\{(x_n,z_n)\}$ which converges to $(x,z)$
with $x_n \in (f_F)^{-1}(ONE_G)$ and   $z_n = f_G \circ f_F(x_n)$.

$\Box$ \vspace{.5cm}

\begin{df}\label{df3.11}  For suitable relations $F : X_1 \to X_2$ and $G : X_2 \to X_3$ define
the \emph{suitable composition} $G \bullet F : X_1 \to X_3$ by
\begin{equation}\label{3.7}
G \bullet F \quad =_{def} \quad \pi_{13}[ \ \overline{(G \otimes F) \cap (X_1 \times ONE_G \times X_3)} \ ].
\end{equation}
\end{df}

The name is justified by the following

\begin{theo}\label{theo3.12} For suitable relations $F : X_1 \to X_2$ and $G : X_2 \to X_3$
the suitable composition $G \bullet F : X_1 \to X_3$ is a suitable
relation with
\begin{equation}\label{3.8}
(f_F)^{-1}(ONE_G) \quad \subset \quad ONE_{G \circ F} \quad \subset \quad ONE_{G \bullet F}. \hspace{1cm}
\end{equation}
$G \bullet F$ is the unique closed subset of $G \circ F$ which is minimal for $\pi_{1G \circ F}$. If $D$ is any
dense subset of $ ONE_{G \circ F}$ then $G \bullet F$ is the closure in $X_1 \times X_3$
of the restriction of $f_{G \circ F}$ to $D$.
\end{theo}

{\bfseries Proof:}The first inclusion of (\ref{3.8}) follows from (\ref{3.6}). On the
other hand, $G \bullet F \subset G \circ F$
and the projection of $G \bullet F$ to $X_1$ is closed and contains the dense set
$(f_F)^{-1}(ONE_G)$.  Hence, for every $x \in X_1$,
$\emptyset \ \not= \ G \bullet F(x) \ \subset \ G \circ F(x)$ which implies the second inclusion of (\ref{3.8}).

Because $ONE_{G \circ F}$ is dense in $X_1$ it follows from Corollary \ref{cor3.2a}
that $G \circ F$ contains a unique minimal subset and that it
is the closure of the map $f_{G \circ F}$ on $ONE_{G \circ F}$ which contains the
map $f_G \circ f_F$ on $(f_F)^{-1}(ONE_G)$. So by Lemma
\ref{lem3.10} it contains $G \bullet F$. $\pi_1$ maps $f_G \circ f_F$ onto the dense set
$(f_F)^{-1}(ONE_G)$ and so
 $G \bullet F$ is mapped onto
$X_1$ by $\pi_1$. It follows that the  minimal set equals $G \bullet F$.
Thus, $\pi_{1G \bullet F}$ is irreducible because it is the restriction of $\pi_{1G \circ F}$ to the minimal subset.

By Theorem \ref{theo3.9}  $f_{G \circ F}$ is weakly open and so by
the Variation of Domain  Theorem \ref{theo2.4} the extension $f_{G \bullet F}$ is weakly open.
By Theorem \ref{theo3.4} $\pi_{2G \bullet F}$ is almost open. It follows that $G \bullet F$ is a suitable relation.

The result for general $D$ follows from Corollary \ref{cor3.3a}.

$\Box$ \vspace {.5cm}

\begin{prop}\label{prop3.12a} Let $F : X_1 \to X_2$ and $G : X_2 \to X_3$ be suitable relations and
let $f : D_1 \to X_2$ and $g : D_2 \to X_3$ be continuous maps which are dense in $F$ and $G$ respectively.

(a) If $f(D_1) \subset D_2$ then the continuous map $g \circ f : D_1 \to X_3$ is dense in $G \bullet F$.

(b) If $D_2$ is a $G_{\delta}$ subset of $X_2$ then $f^{-1}(D_2)$ is a dense subset of $D_1$ and the
continuous map $g \circ f : f^{-1}(D_2) \to X_3$ is dense in $G \bullet F$.
\end{prop}

{\bfseries Proof:} By Corollary \ref{cor3.3a} $D_1$ is a dense subset of $ONE_F$ and $f$ is the restriction of
$f_F$ to $D_1$. Since $f_F$ is weakly open, $f$ is weakly open by the Variation of Domain  Theorem \ref{theo2.4}. Similarly for
$g$. It follows that $f^{-1}(D_2) \subset (f_F)^{-1}(ONE_G)$ and $g \circ f : f^{-1}(D_2) \to X_3$ is a restriction
to $f^{-1}(D_2)$ of $f_G \circ f_F : (f_F)^{-1}(ONE_G) \to X_3$ which is in turn dense in $G \bullet F$. In order to apply
Corollary \ref{cor3.3a} the other way it suffices
to show that $f^{-1}(D_2)$ is dense in $D_1$ and so is dense in $X_1$.

If $f(D_1) \subset D_2$ then $f^{-1}(D_2) = D_1$ and so is dense.  This proves (a).

If $D_2$ is a dense $G_{\delta}$ subset of $X_2$ then  Proposition \ref{prop2.2a} implies
$f^{-1}(D_2)$ is dense in $D_1$ because $f$ is weakly open. This proves (b).

$\Box$ \vspace{.5cm}

\begin{theo}\label{theo3.13} (a) If $F : X_1 \to X_2$ is a suitable relation and
$g : X_2 \to X_3$ is an almost open continuous map
then $g \bullet F = g \circ F$.

(b) If  $G : X_2 \to X_3$ is a suitable relation and $f : X_1 \to X_2$ is an open map then $G \bullet f = G \circ f$.
\end{theo}

{\bfseries Proof:} (a) Since $g$ is a map, $ONE_g = X_2$ and so
$(g \otimes F) \cap (X_1 \times ONE_g \times X_3) = g \otimes F$ which projects
to $g \circ F$.

(b) It suffices to show that $(G \otimes f) \cap (X_1 \times ONE_G \times X_3)$ is dense in $G \otimes f $.

Let $(x,y,z) \in G \otimes f$ and let $U \times V \times W \subset X_1 \times X_2 \times X_3$
be an open set containing $(x,y,z)$.
Since $f(x) = y$ and $f$ is an open map $(V \cap f(U)) \times W$ is an open subset
of $X_2 \times X_3$ which contains $(y,z) \in G$.
Since $G$ is suitable the map $f_G$ is dense in $G$ and so there exists
$ (y_1,z_1) \in (f_G) \cap [(V \cap f(U)) \times W]$. That is,
$y_1 \in ONE_G$ and $(y_1,z_1) \in G$. Since $y_1 \in f(U)$ there exists $x_1 \in U$ with $f(x_1) = y_1$.  Thus,
$(x_1,y_1,z_1) \in (G \otimes f) \cap (X_1 \times ONE_G \times X_3) \cap (U \times V \times W)$.

$\Box$ \vspace{.5cm}

In (b) it is necessary that the map be open. Almost open will not suffice.
Let $X_2 = X_3$ be the unit interval and let $G$ be the
closed relation $F_{01}$ defined by (\ref{1.2a}). Now on
$X_1 = [ - \frac{1}{2},0] \cup [1,\frac{3}{2}]$ define $f :X_1 \to X_2$ by
\begin{equation}\label{3.9}
f(x) = \begin{cases} \ x + \frac{1}{2} \qquad \mbox{for} \quad - \frac{1}{2} \leq x \leq 0, \\
 \ x - \frac{1}{2} \qquad \mbox{for} \quad 1 \leq x \leq \frac{3}{2}. \end{cases}
 \end{equation}
Clearly, $f$ is an irreducible map and $ G \bullet f$ is the irreducible map given by
\begin{equation}\label{3.10}
G \bullet f(x) = \begin{cases} \ \ - x \ \qquad \mbox{for} \quad - \frac{1}{2} \leq x \leq 0, \\
 \  2 - x  \qquad \mbox{for} \quad 1 \leq x \leq \frac{3}{2}. \end{cases}
 \end{equation}
while
 \begin{equation}\label{3.11}
 G \circ f \quad = \quad G \bullet f \cup \{ (0,1), (1,0) \}. \hspace{.5cm}
 \end{equation}
 Note also that $ONE_{G \circ f} = [ - \frac{1}{2},0) \cup (1,\frac{3}{2}]$ is
 a proper subset of $ONE_{G \bullet f} = X_1$.

We will call a relation $F : X_1 \to X_2$ \emph{surjective} when $F(X_1) = X_2$ and $F^{-1}(X_2) = X_1$.  That is,
$F$ and $F^{-1}$ both have full domain.

\begin{prop}\label{prop3.13a} (a) For a suitable relation $F : X_1 \to X_2$ the following are equivalent

\begin{itemize}
\item[(i)]  $F$ is surjective.

\item[(ii)] $\pi_{2F} : F \to X_2$ is a surjective map.

\item[(iii)] $f_F(ONE_F)$ is a dense subset of $X_2$.

\item[(iv)] If $D$ is a dense subset of $ONE_F$ then $f_F(D)$ is a dense subset of $X_2$.
\end{itemize}

(b) If $F : X_1 \to X_2$ and $G : X_2 \to X_3$ are surjective suitable relations then the suitable relation
$G \bullet F : X_1 \to X_3$ is surjective.
\end{prop}

{\bfseries Proof:} Since $F$ is suitable, $\pi_{1F}$ is surjective. So (i) $\Leftrightarrow$ (ii) is obvious.
If $D$ is dense in $ONE_F$ then $f_F(D)$ is dense in $f_F(ONE_F)$ because $f_F$ is continuous.
So (iii) $\Leftrightarrow$ (iv) is obvious.

(iii) $\Rightarrow$ (ii) $F(X_1)$ is closed and contains  $f_F(ONE_F)$.  If the latter is dense then
$F(X_1) = X_2$ which implies (ii).

(ii) $\Rightarrow$ (iii) $f_F$ is dense in $F$ and $\pi_{2F}(f_F) = f_F(ONE_F)$. So $f_F(ONE_F)$ is dense in
$F(X_1)$ in any case.

(b) $D = (f_F)^{-1}(ONE_G)$ is dense in $ONE_F$ and so if $F$ is surjective, $f_F(D) =  ONE_G  \cap f_F(ONE_F)$
is dense in $X_2$ by (a). If $G$ is surjective then $f_G(f_F(D)) = f_{G \bullet F}(D)$ is dense in $X_2$.
This is a subset of $f_{G \bullet F}(ONE_{G \bullet F})$ and so the latter is dense in $X_2$.  By (a) again
$G \bullet F$ is surjective.

$\Box$ \vspace{.5cm}

\begin{theo}\label{theo3.14} If $F : X_1 \to X_2$ , $G : X_2 \to X_3$ and $H : X_3 \to X_4$ are suitable relations then
$(H \bullet G) \bullet F \ = \ H \bullet (G \bullet F) : X_1 \to X_4$.
That is, suitable composition is associative.\end{theo}

{\bfseries Proof:}  Using Proposition \ref{prop3.12a} it is easy to
check that both are the closure of $f_H \circ f_G \circ f_F$ defined on the dense $G_{\delta}$ subset of $X_1$
obtained by pulling back the dense $G_{\delta}$ subset $(f_G)^{-1}(ONE_H)$ of $ONE_G$ via $f_F : ONE_F \to X_2$.

$\Box$ \vspace{.5cm}

It follows that we can define a category whose objects are compact metric spaces
and morphisms are suitable relations. Associativity follows
from Theorem \ref{theo3.14}. By Theorem \ref{theo3.13} \ \
$1_{X_2} \bullet F = 1_{X_2} \circ F = F = F \circ 1_{X_1} = F \bullet 1_{X_1}$.

The definition of isomorphism
in the category says that a suitable relation $F : X_1 \to X_2$ is an isomorphism iff there exists a
suitable relation $G : X_2 \to X_1$ such that $G \bullet F = 1_{X_1}$ and
$F \bullet G = 1_{X_2}$. $G$ is the inverse relation in the category.  When it exists it is unique.

\begin{theo}\label{theo3.15} For a suitable relation $F : X_1 \to X_3$ the following are equivalent.
\begin{itemize}
\item[(i)] $F$ is an isomorphism in the  suitable relations category.

\item[(ii)] $\pi_{2F} : F \to X_2$ is an irreducible map.

\item[(iii)] $F^{-1}$ is a suitable relation.

\item[(iv)] There exist subsets $D_1 \subset X_1$ and $D_2 \subset X_2$ both dense
 and a homeomorphism $f : D_1 \to D_2$ such that $f : D_1 \to X_2$ is dense in $F$.
\end{itemize}

When these conditions hold, the suitable relation $F^{-1}$ is the inverse of $F$ in the suitable relations
 category.  Furthermore,
if $f : D_1 \to D_2$ satisfies the conditions of (iv) then $f^{-1} : D_2 \to X_1$ is dense in $F^{-1}$.
The sets $D_1, D_2$ in (iv) can be chosen to be dense $G_{\delta}$ sets.
\end{theo}

{\bfseries Proof:} (i) $\Rightarrow$ (iii) Let $G$ be the inverse in the category.
$(f_G)^{-1}(ONE_F) \subset ONE_G$ is dense and  $f_F \circ f_G : (f_G)^{-1}(ONE_F) \to X_2$ is dense
in $F \bullet G = 1_{X_2}$. If $y \in (f_G)^{-1}(ONE_F)$ then $f_F(f_G(y)) = y$. Hence,
$f_F(ONE_F)$ contains the dense set $(f_G)^{-1}(ONE_F)$ and is contained in the closed set $F(X_1)$.
Thus, $F(X_1) = X_2$ and so $\pi_{2F}$ is surjective, or, equivalently, $\pi_{1F^{-1}}$ is surjective

 The restriction of
$f_F$ to $(f_F)^{-1}(ONE_G)$ is dense in $F$. Hence, its inverse is dense in $F^{-1}$.

 $f_G \circ f_F : (f_F)^{-1}(ONE_G) \to X_1$ is dense
in $G \bullet F = 1_{X_1}$. Thus, if $x \in (f_F)^{-1}(ONE_G)$ and $y = f_F(x)$ then $f_G(y) = x$.
This implies that the inverse of $f_F : (f_F)^{-1}(ONE_G) \to X_2$ is contained in the closed relation
$G$. Since this inverse is dense in $F^{-1}$.
It follows that $F^{-1} \subset G$.  But $G$ is suitable and $\pi_{1F^{-1}}$ is surjective. From Lemma \ref{lem3.3b}
it follows that $F^{-1} = G$.

This proves (iii) and also shows that the inverse $G$, when it exists, is equal to $F^{-1}$.

(ii) $\Leftrightarrow$ (iii) Let $T : X_1 \times X_2 \to X_2 \times X_1$ be the homeomorphism which switches
coordinates.
\begin{equation}\label{3.12}
\pi_{1F^{-1}} \quad = \quad \pi_{2F} \circ T \qquad \mbox{and} \qquad \pi_{2F^{-1}} \quad = \quad \pi_{1F} \circ T
\end{equation}
Whenever $F$ is suitable
$\pi_{2F^{-1}}$ is irreducible and so is almost open. Hence, $F^{-1}$ is suitable iff the almost open map
$\pi_{1F^{-1}}$ is irreducible.

(iii) $\Rightarrow$ (iv) For any closed relation $F : X_1 \to X_2$ it is always true that:
\begin{equation}\label{3.13}
(ONE_F \times ONE_{F^{-1}}) \cap F \quad = \quad (f_F)^{-1}(ONE_{F^{-1}}) \times (f_{F^{-1}})^{-1}(ONE_F).
\end{equation}
That is, if $(x,y) \in F$ with $x$ and $ y$ in the domains of $f_F$ and $f_{F^{-1}}$ respectively then
$y = f_F(x)$ and $x = f_{F^{-1}}(y)$.

If $F$ and $F^{-1}$ are suitable then $D_1 = (f_F)^{-1}(ONE_{F^{-1}}) \subset X_1$ and
$D_2 =  (f_{F^{-1}})^{-1}(ONE_F) \subset X_2$ are dense $G_{\delta}$ sets and $f_F$ and $f_{F^{-1}}$ are inverse
homeomorphisms between them.

This proves (iv) and also shows that the  sets $D_1$ and $D_2$ can be chosen to be $G_{\delta}$'s.

(iv) $\Rightarrow$ (i) Since the map $f$ is a dense subset of $F$, we can apply the coordinate switching homeomorphism
and see that $f^{-1}$ is dense in $F^{-1}$. Because $f^{-1} : D_2 \to D_1$ is a homeomorphism, $f^{-1} : D_2 \to X_1$
is a weakly open continuous map by Variation of Domain.   By Proposition \ref{prop3.6} $F^{-1}$ is a suitable
relation because it is the closure of a densely defined, weakly open continuous map.

By Proposition \ref{prop3.12a}(a) $1_{D_1} = f^{-1} \circ f $ is dense in $F^{-1} \bullet F$. As the latter is
closed and $D_1$ is dense in $X_1$ we have $1_{X_1} = \overline{1_{D_1}} = F^{-1} \bullet F$.
Similarly, $F \circ F^{-1} = 1_{X_2}$. Thus, $F$ and $F^{-1}$ are inverse
isomorphisms in the suitable relation category.

$\Box$ \vspace{.5cm}

\begin{cor}\label{cor3.16} (a) An almost open continuous map $f: X_1 \to X_2$
between compact metric spaces is an isomorphism
in the suitable relations category iff it is an irreducible map. In that case
\begin{equation}\label{3.13}
1_{X_2} \quad = \quad f \bullet f^{-1} \quad = \quad f \circ f^{-1} \qquad \mbox{and}
\qquad 1_{X_1} \quad = \quad f^{-1} \bullet f.
\end{equation}

(b) If $F : X_1 \to X_2$ is a suitable relation then the almost open continuous map $\pi_{2F} : F \to X_2$ is
a suitable relation and the irreducible map $\pi_{1F} : F \to X_1$ is a suitable relation isomorphism with
inverse $(\pi_{1F})^{-1}$. Furthermore,
\begin{equation}\label{3.14}
F \quad = \quad \pi_{2F} \circ (\pi_{1F})^{-1} \quad = \quad \pi_{2F} \bullet (\pi_{1F})^{-1}.
\end{equation}
If $F$ is a suitable relation isomorphism this expresses $F$ as a suitable composition of an irreducible map and the
inverse of an irreducible map.
\end{cor}

{\bfseries Proof:} (a) For a continuous map $f : X_1 \to X_2$, $\pi_{1f} : f \to X_1$ is a homeomorphism.
Hence, $f$ is irreducible iff $\pi_{2f} : f \to X_2$ is irreducible.

The second equations in (\ref{3.13}) and (\ref{3.14}) follow from Theorem \ref{theo3.13}(a).

The rest should by now be obvious.

$\Box$ \vspace{.5cm}

{\bfseries Remark:} For a surjective map $f: X_1 \to X_2$ $1_{X_1} \ = \  f^{-1} \circ f$
 iff $f$ is injective
and so is a homeomorphism when it is continuous.  So whenever $f$ is an almost one-to-one surjection
 which is not a homeomorphism we can let $G = f^{-1}$ and $F = f$ to get an example where $G \bullet F$ is a proper
 subset of $G \circ F$.
 \vspace{1cm}

\section{Suitable Relation Dynamics}

When the domain and range of a map $f$ or relation $F$ are the same set $E$ we will say $f$ is a map
on $E$ or $F$ is a relation on $E$.  Thus, $F : E \to E$ is a \emph{relation on $E$} and we can iterate
just as we would with a map on $E$. Define $F^0 = 1_E, F^1 = F$ and, inductively, for positive integers $n$
\begin{equation}\label{4.1}
F^{n+1} \quad =_{def} \quad F  \circ F^n \hspace{3cm}
\end{equation}
The inverse relation $F^{-1}$ is a relation on $E$ and we let
\begin{equation}\label{4.2}
F^{-n} \quad =_{def} \quad (F^{-1})^n \quad = \quad (F^n)^{-1}.
\end{equation}

Because composition is associative we have
\begin{equation}\label{4.3}
F^{m+n} \quad = \quad F^{m} \circ F^n
\end{equation}
provided that the integers $m$ and $n$ have the same sign. With opposite signs the result is not true in
general.  For example, $h : E_1 \to E_2$ is a relation then it is easy to check that $h$ is a map iff
\begin{equation}\label{4.4}
h \circ h^{-1} \quad \subset \quad 1_{E_2} \qquad \mbox{and} \qquad h^{-1} \circ h \quad \supset \quad 1_{X_1}
\end{equation}
and the inclusions are equalities iff $h$ is a bijection.

If $A$ is a subset of $E$ we call $A$ a \emph{$+$ invariant subset} for $F$ when $F(A) \subset A$ and
an \emph{invariant subset} when $F(A) = A$.

Now suppose that $F$ is a suitable relation on a compact metric space $X$. We define the suitable iterates by
$F^{\bullet 0} = 1_E, F^{\bullet 1} = F$ and, inductively, for positive integers $n$
\begin{equation}\label{4.5}
F^{\bullet n+1} \quad =_{def} \quad F \bullet F^{\bullet n}
\end{equation}
Since suitable composition is associative we have
\begin{equation}\label{4.6}
F^{\bullet m+n} \quad = \quad F^{\bullet m} \bullet F^{\bullet n}
\end{equation}
for all nonnegative integers $m$ and $n$. By Theorem \ref{theo3.15}
$F^{-1}$ is a suitable relation iff $F$ is an isomorphism in the suitable relations category. In that case, with
\begin{equation}\label{4.7}
F^{\bullet -n} \quad =_{def} \quad (F^{-1})^{\bullet n}\quad = \quad (F^{\bullet n})^{-1},
\end{equation}
equation (\ref{4.6}) holds for all integers $m$ and $n$.

Recall that the continuous, weakly open map $f_F : ONE_F \to X$ is
the restriction of the relation $F$ to the dense $G_{\delta}$
set of points $x$  at which $F(x)$ is a singleton.  With $D_0 = X$ define, inductively, for positive integers
$n$
\begin{equation}\label{4.8}
D_{n} \quad = \quad (f_F)^{-1}(D_{n-1}) \quad = \quad ((f_F)^{n})^{-1}(X).
\end{equation}
That is, $D_n$ is the domain of the iterated relation $(f_F)^n$. Since
$f_F$ is weakly open with a dense domain, each $D_n$ is a dense $G_{\delta}$ subset of $X$ by Proposition
\ref{prop2.2a}. With $n \in {\mathbb N} = \{ 0,1,... \}$ the first equation implies, by induction, that
the sequence $\{ D_n \}$ is monotone decreasing and  we intersect to define
\begin{equation}\label{4.9}
D^{ +}_F \quad =_{def} \quad \bigcap_{n \in {\mathbb N}} \ ((f_F)^{n})^{-1}(X).
\end{equation}

\begin{theo}\label{theo4.1} Let $F$ be a suitable relation on $X$.
$D^{ +}_F \ = \ \bigcap_{n} \ ((f_F)^{n})^{-1}(X)$ is a dense $G_{\delta}$ subset of $X$.
 It is the set of points $x \in ONE_F$ such that  $\{ (f_F)^n(x) \in ONE_F \}$ is defined as an infinite sequence of
 points.   The map $f_F$, or equivalently the relation $F$, restricts to define a map which we denote
 $t_F : D^{+}_F \to D^{+}_F$.

 $F$ is a surjective relation iff $t_F(D^+_F)$ is a dense subset of $D^+_F$ or, equivalently, a dense subset of $X$.

 Let $D \subset X$ and $t : D \to D$ be a continuous map.  The following are equivalent:
 \begin{itemize}
 \item[(i)] Regarded as a map from $D$ to $X$, $t$ is dense in $F$.
 \item[(ii)] $D \subset D^+_F$ dense in $X$ and $t$ is the restriction of $t_F$ to the  set $D$.
 \end{itemize}
 When these conditions hold, $D$ is $+$invariant for $t_F$ and $t^n$ is dense in $F^{\bullet n}$ for
 every positive integer $n$.  In particular, $t_F^n$ is dense in $F^{\bullet n}$ for
 every positive integer $n$.
 \end{theo}

 {\bfseries Proof:} The $G_{\delta}$ set $D^+_F$ is dense in $X$ by the Baire Category Theorem. It is clear
 that if $x \in D^+_F$ then $f_F(x) \in D^+_F$ and so $f_F$ restricts to define the  map $t_F$
 from $D^+_F$ to itself.

Proposition \ref{prop3.13a} (a) implies that $F$ is surjective iff $t_F(D^+_F) = f_F(D^+_F)$ is dense.

By Corollary \ref{cor3.3a} (i) is equivalent to (ii) since $\pi_{1F}$ is irreducible and hence surjective.
Then $t^n$ is dense in $F^{\bullet n}$ for every positive integer $n$ by induction using Proposition \ref{prop3.12a}(a).

$\Box$ \vspace{.5cm}

Thus, we see that the suitable iterations $F^{\bullet n}$ provide the natural closure of the dynamics
given by the iterations $t^n_F$
of the continuous map $t_F$ on the Polish space $D^+_F \subset ONE_F \subset X$.

Now we consider continuous maps between closed relations.

Let $G$ be a closed relation on a compact metric space $Y$ and $H : Y \to X$ be a closed relation between compact
metric spaces.  Define $H \times H : Y \times Y \to X \times X$ by $\{(x_1,x_2,x_3,x_4):(x_1,x_3), (x_2,x_4) \in H \}$.
The  closed relation $H \times H$ is just the set product of the two relations with the second and third
coordinates switched. The image set $(H \times H)(G) \subset X \times X$ is a closed relation on $X$ and it is
easy to check that the image can be written as a composition of closed relations:
\begin{equation}\label{4.10}
(H \times H)(G) \quad = \quad H \circ G \circ H^{-1}. \hspace{2cm}
\end{equation}
If $F$ is a closed relation on $X$ then we can apply this to $H^{-1}$ to get
\begin{equation}\label{4.10a}
(H \times H)^{-1}(F) \quad = \quad H^{-1} \circ F \circ H. \hspace{2cm}
\end{equation}

If $G$ and $F$ are a closed relations on $Y$ and $X$ and $h : Y \to X$ is a continuous map, then
we will say that \emph{$h$ maps $G$ to $F$} when $(h \times h)(G) \subset F$.

\begin{prop}\label{prop4.2}  Let $G, F$ be closed relations on the compact metric spaces $Y$ and $X$ respectively. Let
$h : Y \to X$ be a continuous map. $h$ maps $G$ to $F$ iff
\begin{equation}\label{4.11}
h \circ G \quad \subset \quad F \circ h.
\end{equation}

If $h$ maps $G$ to $F$ then it maps $G^n$ to $F^n$ for every integer $n$.  In particular, $h$ maps $G^{-1}$ to
$F^{-1}$.

If $h$ maps $G$ to $F$ and  $h$ is surjective, then $F$ has full domain (or is surjective) when
$G$ has full domain (resp. is surjective).

Assume $h$ and $\pi_{1G}$ are surjective and that $\pi_{1F}$ is irreducible.
If $h$ maps $G$ to $F$ then $ (h \times h)(G)  = F$.
\end{prop}

{\bfseries Proof:}  If $h \circ G \circ h^{-1} \subset F$ then by (\ref{4.4})
\begin{equation}\label{4.12}
h \circ G \quad \subset \quad h \circ G \circ h^{-1} \circ h \quad \subset \quad F \circ h,
\end{equation}
while if $h \circ G \subset F \circ h$ then by (\ref{4.4}) again
\begin{equation}\label{4.13}
h \circ G \circ h^{-1}  \quad \subset \quad F \circ h \circ h^{-1} \quad \subset \quad F.
\end{equation}

If $h \circ G \subset F \circ h$ then inductively, we have
\begin{equation}\label{4.14a}
h \circ G^n \ = \ h \circ G \circ G^{n-1} \ \subset \ F \circ h \circ G^{n-1 } \
\subset \ F \circ F^{n-1}\circ h \ = \ F^n \circ h.
\end{equation}

Alternatively beginning with $h \circ G \circ h^{-1} \subset F $, we can use (\ref{4.4}) and induction
\begin{equation}\label{4.14}
\begin{split}
h \circ G^n \circ h^{-1}\ = \ h \circ G \circ G^{n-1} \circ h^{-1} \ \subset \hspace{1cm}\\
h \circ G \circ h^{-1} \circ h \circ G^{n-1} \circ h^{-1} \ \subset \ F \circ F^{n-1} \  = \ F^n .
\end{split}
\end{equation}
So $h$ maps $G^n$ to $F^n$ for all positive integers $n$.
Also
\begin{equation}\label{4.15}
h \circ G^{-1} \circ h^{-1} \ = \ (h \circ G \circ h^{-1})^{-1} \ \subset \ F^{-1}.
\end{equation}
That is, $h$ maps $G^{-1}$ to $F^{-1}$. Hence, $h$ maps $(G^{-1})^n$ to $(F^{-1})^n$ for all positive integers $n$.

If $h$ is surjective and $G(Y) = Y$ then $X = h(G(h^{-1}(X))) \subset F(X)$ and if $G^{-1}(Y) = Y$ then
$X = h(G^{-1}(h^{-1}(X))) \subset F^{-1}(X)$. Hence $F$ is surjective when $h$ and $G$ are and $F$ has full domain
when $h$ is surjective and $G$ has full domain.

If $h$ and $\pi_{1G}$ are surjective then $(h \circ G \circ h^{-1})^{-1}(X) = X$ and so with $\tilde F =
h \circ G \circ h^{-1}$ we have $\pi_{1 \tilde F}$ is surjective. $\tilde F \subset F$ since $h$ maps $G$ to $F$.
If $\pi_{1F}$ is irreducible then $\tilde F = F$ by Lemma \ref{lem3.3b} $\tilde F = F$.

$\Box$ \vspace{.5cm}

\begin{theo}\label{theo4.3}   Let $G, F$ be suitable relations on
the compact metric spaces $Y$ and $X$ respectively. Let
$h : Y \to X$ be an almost open, continuous map. $h$ maps $G$ to $F$ iff
\begin{equation}\label{4.16}
h \bullet G \quad = \quad F \bullet h.
\end{equation}

If $h$ maps $G$ to $F$ then it maps $G^{\bullet n}$ to $F^{\bullet n}$ for every nonnegative integer $n$.

If $h$ maps $G$ to $F$ then $D^+_G \cap h^{-1}(D^+_F)$ is a $G_{\delta}$ dense in $Y$ and $+$invariant for
$t_G$. On this set, $h \circ t_G = t_F \circ h$.
\end{theo}

{\bfseries Proof:} By Theorem \ref{theo3.13} (a) $h \circ G = h \bullet G$. Hence, given (\ref{4.16}) we have
\begin{equation}\label{4.17}
h \circ G  \ = \  h \bullet G \  = \  F \bullet h \  \subset \ F \circ h.
\end{equation}

On the other hand, if $h \circ G \subset F \circ h$ then $h \bullet G $ and $  F \bullet h$ are both suitable
relations from $Y$ to $X$ which are contained in $F \circ h$. So each is minimal for the surjection
$\pi_{1 F \circ h}$. By Theorem \ref{theo3.12}, $F \bullet h$ is the unique subset minimal for $\pi_{1 F \circ h}$
and hence $h \bullet G  =  F \bullet h$.

From (\ref{4.16}) we obtain, inductively, for all positive integers $n$
\begin{equation}\label{4.18}
h \bullet G^{\bullet n} \ = \ h \bullet G \bullet G^{\bullet n-1} \ = \ F \bullet h \bullet
G^{\bullet n-1 } \ = \ F \bullet F^{\bullet n-1} \bullet h \ = \ F^{\bullet n} \bullet h.
\end{equation}

Since $h$ is almost open, the $G_{\delta}$ set $D^+_G \cap h^{-1}(D^+_F)$ is dense in $Y$. If
$x \in D^+_G \cap h^{-1}(D^+_F)$ then $t_G(x) \in D^+_G$ since the latter is $t_G$ $+$invariant, and
$h \circ G(x) \subset F \circ h(x)$ which consists of the single point $t_F(h(x))$ because $x \in h^{-1}(D^+_F)$.
Thus, $t_G(x) \in D^+_G \cap h^{-1}(D^+_F)$ and $h(t_G(x)) = t_F(h(x))$.

$\Box$ \vspace{.5cm}

For a closed relation $F$ on a compact metric space $X$ we define the \emph{sample path space} $X^+_F$ to be
\begin{equation}\label{4.19}
X^+_F \quad =_{def} \quad \{ z \in X^{\mathbb N} : (z_n, z_{n+1}) \in F
\quad \mbox{for all } \quad n \in {\mathbb N} \ \}.
\end{equation}
If the domain of $F$ is $X$ then we clearly have:
\begin{equation}\label{4.20}
\pi_0(X^+_F) \ = \ X \qquad \mbox{and} \qquad \pi_0 \bigtriangleup \pi_n(X^+_F) \ = \ F^n,
\end{equation}
for $n = 1,2,...$, where $\pi_0 \bigtriangleup \pi_n$ is defined by $z \mapsto (z_0,z_n)$.

On $X^{\mathbb N}$ we define the \emph{shift map} $\sigma$ by
\begin{equation}\label{4.21}
\sigma(z)_n \quad = \quad z_{n+1} \qquad \mbox{for} \quad n \in {\mathbb N}.
\end{equation}

It is clear that $X^+_F$ is a closed subset which is $+$invariant for the shift. So $\sigma$ restricts to define
a continuous map $\sigma_F$ on $X^+_F$.  The continuous surjection $\pi_0 : X^+_F \to X$ maps $\sigma_F$ on $X^+_F$ to $F$ by
(\ref{4.20}). It is easy to check that the map $\sigma_F$ is surjective on $X^+_F$ iff the relation $F$ is surjective.

The construction is functorial.  If $h : Y \to X$ is a continuous map then the continuous map $h^{\mathbb N} :
Y^{\mathbb N} \to X^{\mathbb N}$ is defined by using $h$ on each coordinate. If $G$ is a closed relation on $Y$
and $h$ maps $G$ to $F$ then it is clear that $h^{\mathbb N}(Y^+_G) \subset X^+_F$ and $h^{\mathbb N}$ maps
$\sigma_G$ to $\sigma_F$.

If $F$ is a suitable relation then we define the \emph{orbit map} $o^+_F : D^+_F \to X^+_F$ by
\begin{equation}\label{4.22}
o^+_F(x)_n \quad = \quad (t_F)^n(x)  \qquad \mbox{for} \quad n \in {\mathbb N}.
\end{equation}

\begin{theo}\label{theo4.4} Let $F$ be a suitable relation on $X$.

(a) The map $o^+_F$ is a homeomorphism from $D^+_F$ onto $X^+_F \cap (\pi_0)^{-1}(D^+_F)$.  In particular,
for $x \in D^+_F \ o_F^+(x)$ is the unique point of $X^+_F$ with $x = o^+_F(x)_0$.

(b) The closed set $S^+_F \subset X^+_F$ defined by
\begin{equation}\label{4.23}
S^+_F \quad = \quad \overline{o^+_F(D^+_F)} \quad = \quad \overline{X^+_F \cap (\pi_0)^{-1}(D^+_F)}
\end{equation}
is the unique subset of $X^+_F$ which is minimal for the restriction of $\pi_0$ to $X^+_F$. For
$n \in {\mathbb N}$ let $\rho_n : S^+_F \to X$ denote the restriction of $\pi_n$ to $S^+_F$.
The map $\rho_0 : S^+_F \to X$ is irreducible and for $n = 1,2,...$
\begin{equation}\label{4.24}
\rho_0 \bigtriangleup \rho_n(S^+_F) \ = \ F^{\bullet n}.
\end{equation}

(c) $S^+_F$ is a closed $+$invariant subset of $X^+_F$. Let $s_F$ denote the restriction of $\sigma_F$ to $S^+_F$.
The map $s_F$ is almost open and $\rho_0$ maps $s_F$ to $F$.

(d) The relation $F$ is surjective iff $s_F$ is a surjective map on $S^+_F$.
\end{theo}

{\bfseries Proof:}  (a) If $x \in D^+_F$ then $(f_F)^n(x) \in ONE_F$ for all $n \in {\mathbb N}$ and so it
is clear that $ o_F^+(x)$ is the unique point of $X^+_F$ with $x = o^+_F(x)_0$. By definition of the
product topology, continuity of $f_F$ implies continuity of $o^+_F$. The continuous map $\pi_0$ on
$X^+_F \cap (\pi_0)^{-1}(D^+_F)$ is the inverse map and so $o^+_F$ is a homeomorphism onto its image.

(b) By part (a) the image $o^+_F(D^+_F)$ is contained in $IN_{\pi_0|X^+_F}$. Since $D^+_F $ is dense in $X$
it is dense in $\pi_0(IN_{\pi_0|X^+_F})$ and the latter is dense in $X$. By Proposition \ref{prop2.19} the
closure of $IN_{\pi_0|X^+_F}$ is the unique subset of $X^+_F$ which is minimal for $\pi_0$ on $X^+_F$.
Since $(\pi_0|X^+_F)^{-1}$ is continuous on $\pi_0|X^+_F(IN_{\pi_0|X^+_F})$ and this continuous map restricts
to $o^+_F$ on the dense set $D^+_F$ it follows that the minimal set is the closure of the image of $o^+_F$.
The restriction $\rho_0$ is irreducible because $S^+_F$ is minimal.

(c) The set $o^+_F(D^+_F)$ is clearly $+$invariant for the shift map and so its closure is as well. The homeomorphism
$\rho_0 : o^+_F(D^+_F) \to D^+_F$ maps the restriction of the shift to the weakly open map $t_F$ on $D^+_F$. Hence,
the restriction of the shift to $o^+_F(D^+_F)$ is weakly open as well.  By the Variation of Domain Theorem
\ref{theo2.4} $s_F$ is weakly open on $S^+_F$. By Theorem \ref{theo2.3} it is almost open. Since $s_F$ is a continuous
map and  $o^+_F(D^+_F)$ is dense in $S^+_F$ it follows that the restriction of $s_F|o^+_F(D^+_F)$ is dense in
$s_F$. Clearly, $\rho_F$ maps  $s_F|o^+_F(D^+_F)$ to $t_F$ on $D^+_F$. Since
$t_F$ is a subset of the closed set $F$ and $s_F|o^+_F(D^+_F)$ is dense in $s_F$ it follows that
$\rho_0 \times \rho_0(s^+_F) \subset F$.

(d) If $F$ is surjective then by Theorem \ref{theo4.1} the set $t_F(D^+_F)$ is dense in $D^+_F$ and hence
$o^+_F(t_F(D^+_F)) = s_F(o^+_F(D^+_F))$ is dense in $o^+_F(D^+_F)$ which is dense in $S^+_F$. Since $s_F(S^+_F)$ is
closed and contains $s_F(o^+_F(D^+_F))$ it equals $S^+_F$.

If $s_F$ is a surjective map on $S^+_F$ then it is a surjective relation and $\rho_0$ is surjective and maps
$s_F$ to $F$.  Hence, $F$ is a surjective relation by  Proposition \ref{prop4.2}.

$\Box$ \vspace{.5cm}

\begin{theo}\label{theo4.5} Let $G$ and $F$ be suitable relations on $Y$ and $X$, respectively and let
$h : Y \to X$ be an almost open continuous map which maps $G$ to $F$. $h^{\mathbb N} : Y^+_G \to X^+_F$
restricts to an almost open map from $S^+_G \to S^+_F$. The restriction maps $s_G$ to $s_F$.\end{theo}

{\bfseries Proof:}  By Theorem \ref{theo4.3} the restriction $h : D^+_G \cap h^{-1}(D^+_F) \to D^+_F$ maps
the restriction of $t_G$ to $t_F$. These are dense sets in $Y$ and $X$ and $h$ is almost open and so weakly open.
By the Variation of Domain Theorem the restriction $h : D^+_G \cap h^{-1}(D^+_F) \to D^+_F$ is weakly open.
Since $o^+_G$ and $o^+_F$ are homeomorphisms, the restriction
$h^{\mathbb N} : o^+_G(D^+_G \cap h^{-1}(D^+_F) ) \to o^+_F(D^+_F)$ is weakly open and maps the restriction
of $s_G$ to the restriction of $s_F$. By Variation of Domain again $h^{\mathbb N} : S^+_G \to S^+_F$ is
weakly open and hence almost open by Theorem \ref{theo2.3}. Since $o^+_G(D^+_G \cap h^{-1}(D^+_F) )$ is dense
in $S^+_G$ and $s_F$ is closed it follows that $h^{\mathbb N}$ maps $s_G$ to $s_F$.

$\Box$ \vspace{.5cm}

\begin{cor}\label{cor4.6} Let $F$ be a suitable relation on $X$ and $g$ be an almost open, continuous map
on a compact metric space $Y$.  If $h : Y \to X$ is an almost open continuous map which maps $g$ to $F$ then
there exists a unique continuous map $h_1 : Y \to S^+_F$ which maps $g$ to $s_F$. The map $h_1$ is almost open
and satisfies $\rho_0 \circ h_1 = h$.\end{cor}

{\bfseries Proof:} $h^{\mathbb N} : S^+_g \to S^+_F$ is an almost open map taking $s_g$ to $s_F$ and $\rho_0 : S^+_g \to Y$ maps $s_g$ to $g$.  Because $g$ is a map $D^+_g = ONE_g = Y$ and the orbit map $o^+_g : Y \to S^+_g$ is a
homeomorphism inverse to $\rho_0$ on $S^+_g$. Hence, $h_1 = h^{\mathbb N} \circ o^+_g$ is the required almost open
continuous map. Uniqueness is obvious on the dense set $h^{-1}(D^+_F)$ and so by continuity uniqueness follows on all
of $Y$.

$\Box$ \vspace{.5cm}

There is a category theory description of what is happening here.  The \emph{suitable dynamics category} has as objects
suitable relations on compact metrizable spaces and has as morphisms almost open maps between them. The \emph{almost open dynamics
category} is the full subcategory whose objects are almost open maps on compact metrizable spaces. The association
from $F$ on $X$ to $s_F$ on $S^+_F$ is a functor from the former category to the latter which is adjoint to the inclusion
functor of the latter into the former.

Let $F$ be a suitable relation on $X$ and $g$ be a continuous map on a compact metric space $Y$.  We will say
that a continuous map $h : Y \to X$ \emph{resolves the discontinuities of $F$ via $g$} if $h$ maps $g$ to $F$
and there exists a dense  $D \subset Y$ which is $+$invariant for $g$ such that $h$ restricts to a
homeomorphism $h : D \to h(D)$ where $h(D)$ is a dense subset of $ONE_F$.

\begin{cor}\label{cor4.7} Let $F$ be a suitable relation on $X$. If $h : Y \to X$ resolves the discontinuities
of $F$ via the continuous map $g$ on $Y$ then $h$ is an irreducible map with $h(D) \subset D^+_F$. $h$ maps
$g$ on $D$ to $t_F$ on $h(D)$. The map $g$ is almost open.
Furthermore, there exists a unique continuous map $h_1 : Y \to S^+_F$ which maps $g$ to $s_F$. The map $h_1$ is irreducible and satisfies $\rho_0 \circ h_1 = h$.
\end{cor}

{\bfseries Proof:}  Because $h : D \to h(D)$ is a homeomorphism both it and its inverse are weakly open. By the usual
Variation of Domain argument, density of $D$ implies that $h$ is almost open and that it is the closure of the
restriction of $h$ to $D$.  Since $h(D)$ is dense, $h$ is surjective and the relation $h^{-1}$ is the closure
of $h^{-1} : h(D) \to D$. By Proposition \ref{prop3.6} the inverse $h^{-1}$ as well as $h$ are suitable relations.
By Corollary \ref{cor3.16} $h$ is irreducible.

Clearly, $h$ maps the orbit sequence of a point $y \in D$ to the orbit sequence of $h(y)$.  That is,
$h(g^n(y)) = (f_F)^n(h(y))$ for every positive integer $n$ and so $h(y) \in D^+_F$ and $h$ maps $g$ on $D$
to $t_F$ on $D^+_F$.

By the Variation of Domain Theorem again, the restriction of $t_F$ to $h(D)$ is weakly open and so composing
with the homeomorphism $h|D$ and its inverse we see that $g$ is weakly open on $D$.  As usual, $g$ is almost open
on $Y$.

The existence of $h_1$ now follows from Corollary \ref{cor4.6}.  It is irreducible by Proposition \ref{prop2.18}
$h_1$ is irreducible because $h$ is.

$\Box$ \vspace{.5cm}

Thus, $\rho_0$ resolves the discontinuities of $F$ via $s_F$ and Corollary \ref{cor3.16} says that every other resolution factors through this one.

In resolving the discontinuities the space $X$ is  replaced by $S^+_F$ which is usually more topologically complicated
as the points of $X$ are ``split'' the relation $F$ in $S^+_F$. However, when $F$ is simple, the space  $Y$ can be simple as well.
Consider from the Introduction the suitable relation $F_{01}$ on the unit interval. Define
$Y = ([0,\frac{1}{2}] \times \{ 0 \}) \ \cup \ ([\frac{1}{2},1] \times \{ 1 \})$. Define the homeomorphism $g$ on $Y$
by
\begin{equation}\label{4.24a}
g(x,i) \quad = \quad (i + \frac{1}{2} - x, i) \qquad \mbox{for} \quad i = 0,1.
\end{equation}
That is, $g$ flips each interval about its midpoint. The projection $\pi_1 : Y \to [0.1]$ maps $g$ into $F$ and resolves the
 discontinuities of the quasi-continuous maps $f_0$ and $f_1$.

We can sharpen the above results when $F$ is an isomorphism in the suitable relations category. If $F$ is a suitable relations isomorphism on $X$ then
(\ref{3.13}) says
\begin{equation}\label{4.25}
(ONE_F \times ONE_{F^{-1}}) \cap F \quad = \quad (f_F)^{-1}(ONE_{F^{-1}}) \times (f_{F^{-1}})^{-1}(ONE_F).
\end{equation}
That is, if $(x,y) \in F$ with $x$ and $ y$ in the domains of $f_F$ and $f_{F^{-1}}$ respectively then
$y = f_F(x)$ and $x = f_{F^{-1}}(y)$ and these sets are dense.

The definition of the dense invariant set $D_F$ on which $f_F$ is a homeomorphism will require a bit
of delicacy.  With $\tilde D_0 = X$ define, inductively, for positive integers
$n$ the monotone decreasing sequence of sets
\begin{equation}\label{4.26}
\tilde D_{n} \quad = \quad (f_F)^{-1}(\tilde D_{n-1}) \cap (f_{F^{-1}})^{-1}(\tilde D_{n-1})
\end{equation}
 Since $f_F$ and $f_{F^{-1}}$ are weakly open with a dense domain, each $\tilde D_n$ is a dense $G_{\delta}$ subset of $X$ by Proposition
\ref{prop2.2a}. With $n \in {\mathbb N} = \{ 0,1,... \}$ we intersect to define
\begin{equation}\label{4.27}
D_F \quad =_{def} \quad \bigcap_{n \in {\mathbb N}} \ \tilde D_{n} \quad \subset \quad D^+_F \cap D^+_{F^{-1}}.
\end{equation}
The final inclusion is obvious but it might be proper.  For example, if $F = f$ is an irreducible map then
$D^+_F = X$ but $D^+_{f^{-1}}$ need not be $+$invariant for $f$. A point $x$ is in $D_F$ when the sequences
$\{ (f_F)^n(x) \}$ and $\{ (f_{F^{-1}})^n(x) \}$ remain in both the domain of $f_F$ and $f_{F^{-1}}$ for every
positive integer $n$. The following is clear and we omit the proof.

\begin{theo}\label{theo4.8} Let $F$ be a suitable relations isomorphism on $X$.
$D_F $ is a dense $G_{\delta}$ subset of $X$.
 The map $t_F$, or equivalently the relation $F$, restricts to define a homeomorphism on $D_F$ whose inverse
 is the restriction of $t_{F^{-1}}$.
 \end{theo}
 \vspace{.5cm}

For a closed relation $F$ on a compact metric space $X$ we define the \emph{sample path space} $X_F$ when $F$ is
surjective, i.e. $F(X) = F^{-1}(X) = X$.
\begin{equation}\label{4.28}
X_F \quad =_{def} \quad \{ z \in X^{\mathbb Z} : (z_n, z_{n+1}) \in F
\quad \mbox{for all } \quad n \in {\mathbb Z} \ \}.
\end{equation}
Since $F$ is surjective we have:
\begin{equation}\label{4.29}
\pi_0(X^+_F) \ = \ X \qquad \mbox{and} \qquad \pi_0 \bigtriangleup \pi_n(X^+_F) \ = \ F^n,
\end{equation}
for $n  \in {\mathbb Z}$, where $\pi_0 \bigtriangleup \pi_n$ is defined by $z \mapsto (z_0,z_n)$.

On $X^{\mathbb Z}$ we define the \emph{shift homeomorphism} $\sigma$ by
\begin{equation}\label{4.30}
\sigma(z)_n \quad = \quad z_{n+1} \qquad \mbox{for} \quad n \in {\mathbb Z}.
\end{equation}

It is clear that $X_F$ is a closed subset which is invariant for the shift. So $\sigma$ restricts to define
a homeomorphism $\sigma_F$ on $X_F$.  The continuous surjection $\pi_0 : X_F \to X$ maps $\sigma_F$ on $X^+_F$ to $F$ and its inverse to $F^{-1}$ by
(\ref{4.29}).

If $F$ is a suitable relations isomorphism then we define the \emph{orbit map} $o_F : D_F \to X_F$ by
\begin{equation}\label{4.31}
o_F(x)_n \quad = \quad (t_F)^n(x)  \qquad \mbox{for} \quad n \in {\mathbb Z}.
\end{equation}

\begin{theo}\label{theo4.9} Let $F$ be a suitable relations isomorphism on $X$.

(a) The map $o_F$ is a homeomorphism from $D_F$ onto $X_F \cap (\pi_0)^{-1}(D_F)$.  In particular,
for $x \in D_F \ o_F(x)$ is the unique point of $X_F$ with $x = o_F(x)_0$.

(b) The closed set $S_F \subset X_F$ defined by
\begin{equation}\label{4.32}
S_F \quad = \quad \overline{o_F(D_F)} \quad = \quad \overline{X_F \cap (\pi_0)^{-1}(D_F)}
\end{equation}
is the unique subset of $X_F$ which is minimal for the restriction of $\pi_0$ to $X_F$. For
$n \in {\mathbb Z}$ let $\rho_n : S_F \to X$ denote the restriction of $\pi_n$ to $S_F$.
The map $\rho_0 : S_F \to X$ is irreducible and for $n \in {\mathbb Z}$
\begin{equation}\label{4.33}
\rho_0 \bigtriangleup \rho_n(S_F) \ = \ F^{\bullet n}.
\end{equation}

(c) $S_F$ is a closed invariant subset of $X_F$. Let $s_F$ denote the restriction of $\sigma_F$ to $S_F$.
The map $s_F$ is a homeomorphism and $\rho_0$ maps $s_F$ to $F$.
\end{theo}

{\bfseries Proof:}  As the proof is completely analogous to that of Theorem \ref{theo4.4} we will omit it, leaving the
adjustments to the reader.

$\Box$ \vspace{.5cm}

From Theorem \ref{theo2.22} it follows that every space without isolated point is isomorphic in the suitable relations
category to the Cantor Set. We show that every suitable relation is isomorphic to a continuous almost open map on
the Cantor Set.

\begin{theo}\label{theo4.10} Let $F$ be a suitable relation on $X$. If $X$ has no isolated points then there
exists an almost open continuous map $g$ on the Cantor Set $C$ and an irreducible map $h : C \to X$ which maps
$g$ to $F$.  Furthermore, if $F$ is a suitable relations isomorphism then $g$ can be chosen to be a homeomorphism on
$C$.\end{theo}

{\bfseries Proof:}  By Theorem \ref{theo2.22} there is an irreducible map $h_1 : C \to X$. By Corollary \ref{cor3.16}
$h_1$ is a suitable relations isomorphism and so we can define a suitable relation $F_1$ on $C$ by
$F_1 = (h_1)^{-1} \bullet F \bullet h_1$ and so $h_1 \bullet F_1 = F \bullet h$. By Theorem \ref{theo4.3} $h_1$
maps $F_1$ to $F$.

Apply Theorem \ref{theo4.4} to $F_1$.  We obtain an almost open map $s_{F_1}$ on $S^+_{F_1}$ and a irreducible
map $\rho_0 : S^+_{F_1} \to C$ which maps $s_{F_1}$ to $F_1$. Now $S^+_{F_1}$ is a closed subset of
$C^{\mathbb N}$ and so is a compact, zero-dimensional, metrizable space.  Since $C$ has no isolated points
and  $\rho_0$ is irreducible, Lemma \ref{lem2.21} implies that $S^+_{F_1}$ has no isolated points.  Hence,
there is a homeomorphism $h_3 : C \to S^+_{F_1}$.  Let $g = (h_3)^{-1} \circ s_{F_1} \circ h_3$. Clearly
$g$ is an almost open map on $C$. Furthermore, $h_1 \circ \rho_0 \circ h_3$ is  irreducible and maps $g$ to $F$.

If $F$ is a suitable relations isomorphism then so is $F_1$.  We apply Theorem \ref{theo4.9} to $F_1$ instead of
Theorem \ref{theo4.4}. This replaces the almost open map $s_{F_1}$ on $S^+_{F_1}$ by the homeomorphism
$s_{F_1}$ on $S_{F_1}$.  Then proceed as before to obtain the homeomorphism $g$.

$\Box$ \vspace{1cm}

\section*{References}

\begin{enumerate}

\item[] E. Akin (1993) {\bfseries The general topology of dynamical systems}, Amer. Math. Soc., Providence, R.I.


\item[] A. Crannell and M. Martelli (2000) \emph{Dynamics of quasicontinuous systems}, J. of Difference Equations and Applications
 {\bfseries 6:} 351-361.

\item[] A. Crannell, M. Frantz and M. LeMasurier (2006) \emph{Closed relations and equivalence classes of quasicontinuous functions}
Real Analysis Exchange, {\bfseries 31.2:}409-424.

\end{enumerate}

\end{document}